\documentclass[11pt,a4paper]{article}
\usepackage[utf8]{inputenc}
\usepackage[T1]{fontenc}
\usepackage{amsmath, amssymb, amsthm}
\usepackage{graphicx}
\usepackage{color}
\usepackage{geometry}
\usepackage{booktabs}
\usepackage{authblk}

\usepackage{tikz}
\usetikzlibrary{positioning}

\usepackage[hidelinks]{hyperref} 
\hypersetup{
    colorlinks=true,    
    linkcolor=blue,     
    citecolor=green,    
    urlcolor=green,     
    pdftitle={Vanishing Control Times for Single-Input Bilinear Quantum Systems under Rapid Drift Spectral Growth},
    pdfauthor={Thomas Chambrion (U. Strasbourg) and Cristian Mendico (U. Bourgogne Europe)},
    pdfsubject={ },
    pdfkeywords={QSL, Quantum control, bilinear control, optimal control, Random matrices},
    pdfproducer={LaTeX with hyperref}, 
    pdfcreator={pdfTeX}                
}
\usepackage{orcidlink}

\title{
Vanishing Control Times for Single-Input Bilinear Quantum Systems under Rapid Drift Spectral Growth
}

\author[1]{T. Chambrion \orcidlink{0000-0002-9006-7596}}
\author[2]{C. Mendico \orcidlink{0000-0001-7674-2578}}

\affil[1]{Institut de Recherche Mathématique Avancée
UMR 7501 Université de Strasbourg et CNRS, 
7 rue René Descartes,
67000 Strasbourg, France}
\affil[2]{Université Bourgogne Europe, CNRS, IMB UMR 5584, 21000 Dijon, France}

\date{\today}

\newtheorem{theorem}{Theorem}
\newtheorem{definition}{Definition}
\newtheorem{proposition}[theorem]{Proposition}
\newtheorem{remark}{Remark}
\newtheorem{corollary}[theorem]{Corollary}
\newtheorem{assumption}{Assumption}
\newtheorem{lemma}[theorem]{Lemma}

\begin{document}

\maketitle

\begin{abstract}
This paper investigates the asymptotic behavior, as the Hilbert space dimension tends to 
infinity, of the state-to-state control time of single-input bilinear quantum systems whose 
control operator is drawn from an isotropic (unitarily invariant) random ensemble. Relying 
on deterministic graph-theoretic bounds, namely star cuts of the interaction graph for the 
lower bound and spanning trees for the upper bound, and on integration over the Haar 
measure, we derive high-probability lower and upper bounds on the control time, expressed 
in terms of the empirical spectral variance and of the relative concentration of the drift 
spectrum. 
Whenever the spectrum grows fast enough (for instance polynomially with an exponent larger 
than 3/2, geometrically, or super-exponentially), the upper bound vanishes and the control 
time tends to zero in probability: in these cases, enlarging the Hilbert space 
constructively speeds up control.
\end{abstract}

\setcounter{tocdepth}{2}
\tableofcontents

\section{Introduction}

\subsection{Control time for bilinear quantum systems}

We consider a real sequence of spectral frequencies $\Lambda = (\lambda_k)$, ordered by increasing modulus: $|\lambda_1| \leq |\lambda_2| \leq \cdots$
 from which we construct, for every dimension 
 $d \in \mathbb{N}^*$, a diagonal drift operator 
 $A^{(d)}= i \, \mathrm{diag}(\lambda_1,\ldots,\lambda_d)$ in the Lie algebra 
 $\mathfrak{u}(d)$. For a given $d$, let $B^{(d)} \in \mathfrak{u}(d)$ be a control operator 
 driving the bilinear quantum control system defined on the unit sphere of $\mathbb{C}^d$:

 $$\frac{d}{dt}\psi(t) = \left(A^{(d)} + u(t)B^{(d)}\right)\psi(t).$$

 Let $\tau_{SS}^{(d)}$ denote the state-to-state control time, defined as the supremum, over all pairs $(\psi,\varphi)$ of unit states, of the minimal time
 needed to steer the system from $\psi$ to $\varphi$ under an unbounded, piecewise-continuous real control field $t\mapsto u(t)$.

 The problem we address is the following: can we predict the asymptotic behavior of $\tau_{SS}^{(d)}$ as the dimension $d \to \infty$, relying solely on the intrinsic properties of the spectral sequence $\Lambda$, for a "generic" interaction matrix $B^{(d)}$ (a notion that will be formalized through random matrix theory)?

\subsection{Deterministic Bounds and the Isotropic Random Framework}

Our asymptotic analysis relies on deterministic estimates of the control time $\tau_{SS}^{(d)}$ established in prior work \cite{Harmonique}. These estimates relate the continuous-time control problem to the topological properties of a discrete interaction graph.

\begin{proposition}[Deterministic bounds, Propositions 1 and 5 
in \cite{Harmonique}]\label{PRO_rappel_estimations_tauSS}
Let $B^{(d)} = P D P^\dagger \in \mathfrak{u}(d)$ be a spectral decomposition of the control operator ($D$ is diagonal and $P$ is a unitary matrix), and define the interaction matrix in the eigenbasis of $B^{(d)}$ as $\tilde{A}^{(d)} = P^\dagger A^{(d)} P$. For $j \in \{1,\dots,d\}$, denote by

$$C_j = \sum_{k \neq j} |\tilde{A}^{(d)}_{kj}|^2,$$

the capacity of the star cut isolating the vertex $j$. Let $T$ be any spanning tree of the complete graph on $\{1,\dots,d\}$ equipped with edge weights $|\tilde{A}^{(d)}_{jk}|^{-1}$. Then, for every $j \in \{1,\dots,d\}$, the state-to-state control time satisfies the unconditional lower bound
$ \displaystyle \tau_{SS}^{(d)} \geq \frac{2}{\sqrt{C_j}},$
and, if the spectrum of $D$ is rationally independent, the upper bound
$ \displaystyle \tau_{SS}^{(d)} \leq \frac{\pi}{2} \sum_{(j',k') \in T} \frac{1}{|\tilde{A}^{(d)}_{j'k'}|}$.
\end{proposition}

Since $A^{(d)}$ is skew-Hermitian, $|\tilde{A}^{(d)}_{kj}| = |\tilde{A}^{(d)}_{jk}|$ and the edge weights are symmetric. 
Both bounds are invariant under a translation $\lambda_m \mapsto \lambda_m + c$ of the drift spectrum, since such a translation only modifies the diagonal of $\tilde{A}^{(d)}$. 
Throughout the paper, only star cuts are used for the lower bound: any of the $d$ star cuts can be chosen, and the choice may depend on the realization of $P$.

To predict the typical behavior of $\tau_{SS}^{(d)}$ as $d \to \infty$ independently of pathological worst-case configurations, we evaluate these bounds when the control operator $B^{(d)}$ is drawn from a generic, isotropic distribution. 

By an "isotropic random matrix", we denote any random skew-Hermitian matrix $B \in \mathfrak{u}(d)$ whose probability distribution is unitarily invariant (i.e., lacking any privileged quantization axis) and whose eigenvalues admit a joint continuous density (without atoms). As formally detailed in Section \ref{SEC_Preliminaries}, the Weyl integration formula ensures that the geometry of the interaction and the intensities of the couplings statistically decouple:
\begin{enumerate}
    \item The eigenvector matrix $P \in U(d)$ is distributed according to the uniform Haar measure. This guarantees that the eigenbasis of $B$ points in a random direction within the Hilbert space.
    \item The eigenvalue matrix $D$ admits a joint probability density with respect to the Lebesgue measure.
\end{enumerate}

Since the rational independence condition of Proposition \ref{PRO_rappel_estimations_tauSS} holds almost surely in this framework, the asymptotic behavior of $\tau_{SS}^{(d)}$ is governed by the drift spectrum $\Lambda$ alone; we identify three canonical regimes.

\subsection{Main results}

Throughout the paper, we denote by 
\[
\bar{\lambda}_d = \frac{1}{d}\sum_{k=1}^d \lambda_k \quad \mbox{ and } \quad v_d = \frac{1}{d}\sum_{k=1}^d (\lambda_k-\bar{\lambda}_d)^2 = \frac{1}{d}\sum_{k=1}^d \lambda_k^2 - \bar{\lambda}_d^2,
\]
the empirical spectral mean and the empirical spectral variance of $A^{(d)}$. Note that $v_d$, as well as the bounds of Proposition \ref{PRO_rappel_estimations_tauSS}, are invariant under a translation of the spectrum.

The first regime, called the Lyapunov regime, corresponds to cases where the relative weight of the largest eigenvalue in the empirical spectral variance asymptotically vanishes.
This encompasses spectra $\Lambda$ exhibiting polynomial or logarithmic growth, which represent most of standard quantum mechanical models.

\begin{assumption}[Lyapunov condition]\label{ASS_Lyapunov_condition}
Let $A^{(d)} = i \, \mathrm{diag}(\lambda_1, \dots, \lambda_d)$ be a sequence of deterministic drift matrices with $v_d>0$. 
We say that the eigenvalues satisfy the Lyapunov condition if 
\[\frac{1}{d}\sum_{k=1}^d (\lambda_k-\bar{\lambda}_d)^4 = o(d v_d^2) \quad \mbox{ as } d \to \infty.\]
\end{assumption}
Assumption \ref{ASS_Lyapunov_condition} is invariant under translation of the spectrum. Since $(x-y)^4 \le 8(x^4+y^4)$ and $\bar{\lambda}_d^4 \le \frac{1}{d}\sum_k \lambda_k^4$ (Jensen's inequality), it is implied by the uncentered condition $\frac{1}{d}\sum_{k=1}^d \lambda_k^4 = o(d v_d^2)$.
The Lyapunov condition is for instance satisfied for every spectrum of the form $\lambda_m=  m^a$ with $a >0$ (see Section \ref{SEC_polynomial_spectrum} below).

\begin{theorem}[Time estimates in Lyapunov regime]\label{THE_random_encadrement}

Let $A^{(d)} = i \, \mathrm{diag}(\lambda_1, \dots, \lambda_d)$ be a sequence of deterministic drift matrices with $v_d>0$, satisfying the Lyapunov condition (Assumption \ref{ASS_Lyapunov_condition}).
Then, for any arbitrarily small $\epsilon > 0$, 
 there exist a constant $C_\epsilon > 0$ and an integer $d_0 \in \mathbb{N}$ such that for all $d \ge d_0$:$$\mathbb{P} \left( \frac{2 - \epsilon}{\sqrt{v_d}} \le \tau_{SS}^{(d)} \le C_\epsilon \frac{d^{3/2}}{\sqrt{v_d}} \right) \ge 1 - \epsilon.$$
\end{theorem}

The upper bound in Theorem \ref{THE_random_encadrement} actually holds for every spectrum and every $d$, without the Lyapunov condition (Theorem \ref{TH_Uniform_bound_MST}). Similarly, a weaker lower bound $\tau_{SS}^{(d)} \geq 2\sqrt{\epsilon}/\sqrt{v_d}$, holding with probability at least $1-\epsilon$, is valid for every spectrum and every $d$ (Proposition \ref{PRO_Uniform_lower_bound}). The role of the Lyapunov condition is to provide the improved constant $2$ in the lower bound.

The second regime, called uniformly geometric, corresponds to cases where the ratio of two adjacent eigenvalues of $A^{(d)}$ is bounded away from $1$ and from $+\infty$.  In this case, the relative weight of the largest eigenvalue in the $\ell^2$ norm of the spectrum is bounded away from $0$ and from $1$.   
This includes for instance spectra with an exponential growth $\lambda_d= a^d$, $a>1$.

\begin{assumption}[Uniformly Geometric Spectral Class]\label{ASS_geometric_class}
Let $\Lambda= (\lambda_m)_{m \ge 1}$ be the sequence of eigenvalues of the drift operator, ordered such 
that $|\lambda_1| \le |\lambda_2| \le \dots$. The spectrum belongs to the uniformly geometric 
class $\mathcal{G}(a,b)$ if there exist constants $a$ and $b$ with $1 < a \le b < \infty$, and 
an integer $m_0 \ge 1$, such that for all $m \ge m_0$:
$$a \le \frac{|\lambda_{m+1}|}{|\lambda_m|} \le b.$$
\end{assumption} 

\begin{theorem}[Time Estimates for the Uniformly Geometric Class]\label{THE_random_encadrement_geometric}
Let $A^{(d)} = i \, \mathrm{diag}(\lambda_1, \dots, \lambda_d)$ be a sequence of deterministic drift matrices, ordered such that $|\lambda_1| \le \dots \le |\lambda_d|$.
Assume the spectrum belongs to the Uniformly Geometric Spectral Class (Assumption \ref{ASS_geometric_class}). Let $\tau_{SS}^{(d)}$ denote the state-to-state control time of the randomized system.

For any arbitrarily small failure probability $\epsilon > 0$, there exist strictly positive constants $c_\epsilon$ and $C_\epsilon$ and an integer $d_0 \in \mathbb{N}$ such that, for all $d \geq d_0$:
\[\mathbb{P} \left( c_\epsilon \frac{\sqrt{d}}{|\lambda_d|} \le \tau_{SS}^{(d)} \le C_\epsilon \frac{d^2}{|\lambda_d|} \right) \ge 1 - \epsilon.\]
\end{theorem}

In the uniformly geometric class, $v_d \asymp \lambda_d^2/d$ (Lemma \ref{LEM_variance_geometric}), so that Theorem \ref{THE_random_encadrement_geometric} can be equivalently written $c_\epsilon d^{-1/2}\sqrt{d/v_d} \le \tau_{SS}^{(d)} \le C_\epsilon d \sqrt{d/v_d}$, which is the same scaling as in Theorem \ref{THE_random_encadrement}.

The last case (spectral dominance) covers the case where the relative weight of the largest eigenvalue tends to 1. This happens for instance with superexponential spectra such as $\lambda_d=d!$. 

\begin{assumption}[Absolute Spectral Dominance Assumption]\label{ASS_spectral_dominance}
Assume the spectrum $(\lambda_1, \dots, \lambda_d)$ is ordered such that $|\lambda_1| \le \dots \le |\lambda_d|$ and satisfies:
\[\lim_{d \to \infty} \frac{\lambda_d^2}{\sum_{m=1}^d \lambda_m^2} = 1, \quad \mbox{ or equivalently } \quad \sum_{m=1}^{d-1} \lambda_m^2 = o(\lambda_d^2).\]
\end{assumption}

The lower bound in the dominance regime requires a quantitative version of Assumption \ref{ASS_spectral_dominance}.

\begin{assumption}[Quantitative Spectral Dominance]\label{ASS_quantitative_dominance}
Assume the spectrum $(\lambda_1, \dots, \lambda_d)$ is ordered such that $|\lambda_1| \le \dots \le |\lambda_d|$ and that there exist $K>0$ and $d_0 \in \mathbb{N}$ such that, for every $d \geq d_0$,
\[ \sum_{m=1}^{d-1} \lambda_m^2 \leq K \frac{\lambda_d^2}{d}.\]
\end{assumption}
Assumption \ref{ASS_quantitative_dominance} clearly implies Assumption \ref{ASS_spectral_dominance}. It is satisfied for instance by $\lambda_d = d!$, for which $\sum_{m<d} \lambda_m^2 \sim ((d-1)!)^2 = \lambda_d^2/d^2$, or by $\lambda_d = e^{d^2}$. It is not satisfied by $\lambda_d = \prod_{k=2}^d \log k$, which nevertheless satisfies Assumption \ref{ASS_spectral_dominance} (see Remark \ref{REM_quantitative_dominance_needed}).

\begin{theorem}[Time Estimates for Absolute Spectral Dominance]\label{THE_random_encadrement_dominance}
Let $A^{(d)} = i \, \mathrm{diag}(\lambda_1, \dots, \lambda_d)$ be a sequence of deterministic drift matrices, ordered such that $|\lambda_1| \le \dots \le |\lambda_d|$. Let $\tau_{SS}^{(d)}$ denote the state-to-state control time of the randomized system.
\begin{enumerate}
\item If the spectrum satisfies the Absolute Spectral Dominance condition (Assumption \ref{ASS_spectral_dominance}), then for any $\epsilon > 0$, there exist $C_\epsilon > 0$ and an integer $d_0 \in \mathbb{N}$ such that, for all $d \geq d_0$:
$$\mathbb{P} \left( \tau_{SS}^{(d)} \le C_\epsilon \frac{d^2}{|\lambda_d|} \right) \ge 1 - \epsilon.$$
In particular, $\tau_{SS}^{(d)}$ tends to $0$ in probability as $d \to \infty$.
\item If the spectrum satisfies the Quantitative Spectral Dominance condition (Assumption \ref{ASS_quantitative_dominance}), then for any $\epsilon > 0$, there exist $c_\epsilon, C_\epsilon > 0$ and an integer $d_0 \in \mathbb{N}$ such that, for all $d \geq d_0$:
$$\mathbb{P} \left( c_\epsilon \frac{d}{|\lambda_d|} \le \tau_{SS}^{(d)} \le C_\epsilon \frac{d^2}{|\lambda_d|} \right) \ge 1 - \epsilon.$$
\end{enumerate}
\end{theorem}

\subsection{Illustrative Example: The Polynomial Spectrum 
}\label{SEC_polynomial_spectrum}

To illustrate the practical implications of the topological bounds, we evaluate the asymptotic envelope of the control time $\tau_{SS}^{(d)}$ for a polynomial spectrum $\lambda_k = k^a$ (assuming first $a > 0$).

We must first evaluate the empirical spectral variance $v_d$ and verify the Lyapunov condition. The asymptotic spectral moments can be computed via Riemann sums:
$$S_p = \sum_{k=1}^d k^{pa} \sim \frac{d^{pa+1}}{pa+1}.$$
The empirical variance $v_d = \frac{S_2}{d} - \left(\frac{S_1}{d}\right)^2$ evaluates to:
$$v_d \sim \frac{d^{2a}}{2a+1} - \frac{d^{2a}}{(a+1)^2} = \left[ \frac{a}{(a+1)\sqrt{2a+1}} \right]^2 d^{2a}.$$

The applicability of Theorem \ref{THE_random_encadrement} relies on the Lyapunov condition (Assumption \ref{ASS_Lyapunov_condition}). By the remark following Assumption \ref{ASS_Lyapunov_condition}, it is enough to check the uncentered ratio:
$$\frac{\frac{1}{d}\sum_{k=1}^d \lambda_k^4}{d \cdot v_d^2} \sim \frac{\frac{d^{4a}}{4a+1}}{d \cdot \left[ \frac{a}{(a+1)\sqrt{2a+1}} \right]^4 d^{4a}} = \mathcal{O}\left(\frac{1}{d}\right).$$
Since this ratio vanishes as $d \to \infty$, the Lyapunov condition is satisfied for any exponent $a > 0$.

Defining the spectrum-specific geometric constant $K_a = \frac{a+1}{a}\sqrt{2a+1}$, the injection of $v_d$ into Theorem \ref{THE_random_encadrement} establishes the asymptotic behavior of the control time. For any arbitrarily small $\epsilon > 0$, there exists a constant $C_\epsilon > 0$ such that, for sufficiently large dimensions $d$, the following bounds hold with probability at least $1-\epsilon$:

$$(2-\epsilon) K_a \frac{1}{d^a} \le \tau_{SS}^{(d)} \le C_\epsilon K_a d^{\frac{3}{2} - a}.$$

This envelope reflects the competition between the polynomial dispersion rate of the spectrum (the exponent $a$) and the effective dimension of the Haar interaction graph (which generates the factor $d^{3/2}$ in the numerator). Depending on $a$, the system falls into distinct physical regimes:

\paragraph{Regime of strong controllability ($a > \frac{3}{2}$).}
The upper bound converges to zero. This is the case for the quadratic rigid rotor ($a=2$). The control time tends to zero as the dimension grows, confirming that expanding the Hilbert space constructively facilitates bilinear inversion, provided the spectrum grows sufficiently fast to overcome the $d^{3/2}$ factor.

\paragraph{Critical regime ($a = \frac{3}{2}$).} 
The theoretical upper bound scaling evaluates to an order of $\mathcal{O}(1)$. However, 
this bound is obtained by evaluating a fixed (path) spanning tree, and is very likely not sharp for the minimal spanning tree. 
If the off-diagonal entries $\tilde{A}_{jk}$ were independent, with the same complex Gaussian marginal distribution of variance $v_d/d$, 
the largest coupling of a given vertex would be of order $\sqrt{v_d \ln d / d}$, and the weight of the minimal spanning tree would be of order 
$d^{3/2}/\sqrt{v_d \ln d}$, that is, a logarithmic refinement of Theorem \ref{TH_Uniform_bound_MST}. While the algebraic 
correlations inherent to Haar-distributed matrices preclude the formal application of 
this independent model, numerical simulations in dimensions $d=1000$ to $5000$ 
strongly suggest that the refinement holds, yielding a conjectured upper bound  
scaling as $\tau_{SS}^{(d)} = \mathcal{O}((\ln d)^{-1/2})$ for $a=3/2$, \cite{Numerical}. The derivation of 
this logarithmic factor remains an open problem.

\paragraph{Regime of weak controllability ($0 < a < \frac{3}{2}$).}
This encompasses the linear spectrum of the harmonic oscillator ($a=1$). The upper bound diverges (e.g., $\tau_{SS}^{(d)} \le \mathcal{O}(d^{1/2})$ for $a=1$). This does not formally imply that the physical system is slow, since the lower cut-based bound still converges to zero ($
\tau_{SS}^{(d)} \gtrsim d^{-a}$). However, it demonstrates that the sequential routing strategy (the Minimal Spanning Tree) mathematically fails to guarantee short-time control in this regime.
Characterizing the true control time in this regime would require going beyond the sequential routing strategy.

\paragraph{Degenerate spectrum ($a = 0$).} 
In this extreme scenario, the spectrum is fully degenerate ($\lambda_k = 1$).  The drift matrix $A$ commutes with any unitary transformation $P$, leaving the conjugated interaction matrix $\tilde{A}$ entirely diagonal. The interaction graph is empty (all transverse couplings vanish), and the bilinear system  loses its controllability.

\paragraph{Regime of absolute speed limit ($a < 0$).} 
For negative exponents, the sequence $(k^a)$ is decreasing and accumulates near zero (the ordering $|\lambda_1| \leq |\lambda_2| \leq \cdots$ is then recovered by a reindexation, see Section \ref{SEC_Random_matrices_Haar}). Since all eigenvalues lie in $(0,1]$, the star cut capacities of Proposition \ref{PRO_rappel_estimations_tauSS} are variances of probability distributions supported in $(0,1]$: $C_j = \sum_m \lambda_m^2 |P_{mj}|^2 - (\sum_m \lambda_m |P_{mj}|^2)^2 \leq 1/4$. Hence, deterministically, $\tau_{SS}^{(d)} \ge 4$ for every $d$ and every control operator.

The probabilistic estimate complements this in the sub-range $-1/2 \le a < 0$. 

For $-1/2 < a < 0$ the Riemann-sum computation above still holds and $v_d \asymp d^{2a} \to 0$; at the borderline $a = -1/2$ the second-moment sum diverges logarithmically, $\sum_{k=1}^d k^{-1} \sim \ln d$, and 
$v_d \sim (\ln d)/d$. The Lyapunov condition is satisfied throughout this range (for
$a \le -1/4$ the fourth-moment sum $\sum k^{4a}$ converges, which only accelerates the $o(1)$), hence Theorem \ref{THE_random_encadrement} gives, in probability,
$$\tau_{SS}^{(d)} \ge \frac{2-\epsilon}{\sqrt{v_d}} \gtrsim \begin{cases} d^{-a}, & -1/2 < a < 0, \\ \sqrt{d/\ln d}, & a = -1/2, \end{cases},$$
which diverges as $d \to \infty$: the control time grows without bound. 

For $a < -1/2$ the Lyapunov condition fails (the Lyapunov ratio tends to the non-zero constant $\zeta(-4a)/\zeta(-2a)^2$) and the spectrum lies outside the three canonical regimes. In all cases $a < 0$, expanding the Hilbert space provides no energetic resource, precluding arbitrarily fast control.

\subsection{Contribution and position with respect to the literature}

The present work bridges graph-theoretic quantum speed limits with random matrix 
theory. Our starting point relies on recent deterministic bounds (Proposition \ref{PRO_rappel_estimations_tauSS}) for the state-to-state 
control time, established through combinatorial structures: cuts of the interaction graph (star cuts in the present paper) 
for the lower bound, and effective spanning trees for the upper bound. Our 
methodology evaluates these topological bounds under random Haar unitary conjugations. 
While the probabilistic machinery employed -- ranging from measure concentration to 
extreme value theory -- relies on established mathematical physics frameworks, its 
main contribution to the quantum control community is the following: we demonstrate 
that in quantum systems, specific spectral growth rates (e.g., quadratic 
dispersion) decrease the required control time as the 
system dimension increases.

Our stochastic evaluation of the deterministic topological bounds relies on combining algebraic integration over compact groups with asymptotic random matrix and extreme-value theories.

The  evaluation of the interaction matrix moments and their universal scaling relies 
on integration over the unitary group. We employ the Weingarten calculus \cite{Weingarten1978, Collins2003} to compute polynomial expectations of the Haar measure, 
complemented by measure concentration phenomena for classical compact groups 
\cite{Meckes_2019} to  establish the asymptotic isotropy characterizing the 
Lyapunov regime (for which the relative weight of the largest eigenvalue in the empirical spectral variance asymptotically vanishes).

In the intermediate uniformly geometric class, where the relative weight of the largest eigenvalue neither tends to 0 nor 1, the spectrum is a Hadamard lacunary 
sequence whose growth ratio is bounded both above and below.
The convergence of weighted sums with lacunary coefficients, which is the mechanism behind the limit laws in our geometric regime,  is a classical topic in the theory of random Fourier series
\cite{Kahane_1986}. Geometric spacing suppresses 
structural correlations, yielding asymptotic independence \cite{Zygmund_1959}.  
From a random matrix perspective, the conjugated drift is a bulk spectrum plus a rank-one 
spike -- the dominant eigenvalue in a Haar-random direction. The theory of finite-rank 
perturbations of large random matrices \cite{benaych-georges2011eigenvalues} characterizes 
when such a spike dominates its bulk; in the present setting the dominance regime is 
precisely that case, in which the entries of $\tilde{A}$ are governed by the spike term 
$\lambda_d\,\bar P_{dj} P_{dk}$.

Finally, in the absolute spectral dominance regime, where the largest eigenvalue tends to concentrate the whole norm of the spectrum, the interaction matrix becomes asymptotically of rank one: 
$|\tilde{A}_{jk}| \approx |\lambda_d| |P_{dj}| |P_{dk}|$. The weakest cut of the network isolates the vertex $j$ minimizing $|P_{dj}|$, 
and its capacity is governed by the minimum of $d$ (almost) independent exponential variables. More generally, extreme values of the entries of Haar matrices 
(Fréchet tails for the inverse couplings) determine the bottlenecks of the network 
\cite{leadbetter1983extremes, embrechts1997modelling}. 

We contrast our methodology, which 
relies on bounding the control time via specific spanning trees, with classical 
probabilistic results on complete random graphs. The weight of the Minimum Spanning 
Tree (MST) converges to a macroscopic constant (such as $\zeta(3)$) for complete graphs with 
independent and identically distributed edge weights whose density is positive at $0$ \cite{frieze1985value}. In our setting, the edge weights $|\tilde{A}_{jk}|^{-1}$ are bounded away from zero with high probability, and their scale depends on $d$: 
even in an idealized model where the entries $\tilde{A}_{jk}$ would be independent with the same marginal distribution, the MST weight would be of order $d^{3/2}/\sqrt{v_d \ln d}$ (see the discussion of the critical regime above). The $\mathcal{O}(d)$ factor appearing in our upper bounds (with respect to the typical edge weight $\sqrt{d/v_d}$) is therefore due to the number of edges of a spanning tree and to the scale of the marginal distribution, and not to the algebraic correlations of Haar unitary matrices.

\subsubsection{Content of the paper}

The remainder of this paper is organized as follows. Section \ref{SEC_Preliminaries} introduces the mathematical framework. Subsection \ref{SEC_Random_matrices_Haar} formalizes the isotropic random framework, in particular the factorization of the law of the control operator via the Weyl integration formula; Subsection \ref{SEC_preliminary_moments} establishes the spectrum-independent finite moments of the conjugated drift matrix; Subsection \ref{SEC_Uniform_bound_MST} derives a universal upper bound on the minimal spanning tree weight, together with a universal (crude) lower bound on the control time; and Subsection \ref{SEC_Gaussian_Approximation} recalls the Gaussian approximation of the entries of Haar matrices (Borel--Poincar\'e lemma and its two-column extension), which underlies the limit laws of the couplings in Section \ref{SEC_Classification} .

Section  \ref{SEC_Lyapunov_normality} establishes the asymptotic normality of the interaction matrix entries under the Lyapunov condition. Section \ref{SEC_Time_estimates}  proves the main state-to-state time estimates in the three spectral regimes (Lyapunov, uniformly geometric, spectral dominance). Section  \ref{SEC_Classification} synthesizes these results: it identifies the limit laws of the couplings in the three regimes and shows that the spectral shape vector, rather than the relative concentration alone, is the finer invariant on which the control-time estimates ultimately depend.

\section{Mathematical Framework and Preliminaries}\label{SEC_Preliminaries}
 Throughout the paper, all random variables take values in $\mathbb{R}^n$ or $\mathbb{C}^n \simeq \mathbb{R}^{2n}$. 

 A complex random variable $Z$ is a standard complex Gaussian variable, denoted $Z \sim \mathcal{CN}(0,1)$, if $Z = X + iY$ where $X$ and $Y$ are independent real Gaussian variables with law $\mathcal{N}(0, 1/2)$; then $\mathbb{E}[Z] = 0$ and $\mathbb{E}[|Z|^2] = 1$. More generally, $Z \sim \mathcal{CN}(0,s^2)$ if $Z/s \sim \mathcal{CN}(0,1)$, and a random vector $G \in \mathbb{C}^d$ is a standard complex Gaussian vector, denoted $G \sim \mathcal{CN}(0, I_d)$, if its coordinates are independent $\mathcal{CN}(0,1)$ variables. The law of such a vector is invariant under the action of the unitary group $U(d)$, and $\langle v, G \rangle \sim \mathcal{CN}(0, \|v\|^2)$ for every deterministic $v \in \mathbb{C}^d$.  If $Z\sim\mathcal{CN}(0,1)$, then $|Z|$ has density $2re^{-r^2}$, in particular $\mathbb{E}[|Z|]=\sqrt{\pi}/2$ and $\mathbb{E}[|Z|^{-1}]=\sqrt{\pi}$.

\subsection{Random skew-Hermitian matrices}\label{SEC_Random_matrices_Haar} 

We are interested in $\tau_{SS}^{(d)}$ for the system when the control generator $B^{(d)}$ is drawn from an isotropic (unitarily invariant) ensemble, i.e., a law on $\mathfrak{u}(d)$ that favors no basis. Concretely, the framework is the following factorization.
 \begin{proposition} \label{PRO_factorisation_mesure_Haar}
 Let $B \in \mathfrak{u}(d)$ be a random skew-Hermitian matrix whose probability distribution admits a density $f(B)$ with respect to the Lebesgue measure on $\mathfrak{u}(d)$. Assume that this distribution is unitarily invariant, meaning $f(U B U^\dagger) = f(B)$ for any unitary transformation $U \in U(d)$.
 Consider a spectral decomposition $B = P D P^\dagger$ (the choice being immaterial, as discussed in the proof), where $D = i \, \mathrm{diag}(\omega_1, \dots, \omega_d)$ contains the real spectrum $\omega \in \mathbb{R}^d$, and $P \in U(d)$ is the unitary matrix of eigenvectors.
 Then, the eigenvector matrix $P$ and the eigenvalue spectrum $\omega$ are independent random variables. Furthermore, the probability measure $d\mu(B) = f(B)dB$ factorizes as: 
 $$d\mu(B) = \tilde{f}(\omega)  d\omega  \otimes dP,$$
 where $\tilde{f}(\omega)$ depends only on the spectrum (the eigenvalues being labelled in a uniformly random order),  and $dP$ is the unique normalized Haar measure on $U(d)$ (projected onto the flag manifold $U(d)/\mathbb{T}^d$ to account for the phase invariance of eigenvectors).
 \end{proposition}
\begin{proof}
We begin with an observation on the non-uniqueness of the spectral decomposition. 
For a given $B$, the factorization $B=PDP^\dagger$ is determined only up to the 
permutation of the pairs $(\omega_j, p_j)$, where $p_j$ is the $j$-th column of $P$, 
and to the multiplication of the columns by unimodular constants; when $B$ has multiple 
eigenvalues, the columns of a common eigenspace may in addition be rotated by an 
arbitrary unitary matrix. None of these choices affects the statement: the permutation 
is immaterial because the formula below treats $\omega$ as a labelled $d$-tuple (the 
labelling being uniformly random), the phases are factored out by taking $P$ in the 
flag manifold $U(d)/\mathbb{T}^d$, and, as shown below, $B$ has simple spectrum almost 
surely, so that the eigenspace ambiguity is void; on the null set where it occurs, any 
spectral decomposition may be chosen.

By the Weyl integration formula, valid for unitarily invariant functions (see 
\cite{Meckes_2019}, Chapter 3), there is a constant $C_d>0$ (depending only on the 
normalization of $dB$) such that, for every unitarily invariant $g\in L^1(\mathfrak{u}
(d))$,
\begin{equation}\label{EQ_weyl_integration}
\int_{\mathfrak{u}(d)} g(B)\,dB = C_d\int_{\mathbb{R}^d} g(i\,\mathrm{diag}(\omega_1,
\dots,\omega_d))\prod_{j<k}(\omega_j-\omega_k)^2\,d\omega.
\end{equation}
(When the inner variable is made explicit, this is $\int_{\mathbb{R}^d}\int_{U(d)/\mathbb{T}^d}g(P\Delta P^\dagger)\,dP\prod_{j<k}(\omega_j-\omega_k)^2\,d\omega$, since unitary invariance gives $g(P\Delta P^\dagger)=g(\Delta)$.)

Set $\tilde f(\omega):=C_d f(i\,\mathrm{diag}(\omega))\prod_{j<k}(\omega_j-\omega_k)^2$. Since $f$ is a probability density, it belongs to $L^1(\mathfrak{u}(d))$ and, by assumption, is unitarily invariant; applying \eqref{EQ_weyl_integration} with $g=f$ yields $\int_{\mathbb{R}^d}\tilde f(\omega)\,d\omega=\int_{\mathfrak{u}(d)}f(B)\,dB=1$, so that $\tilde f\,d\omega$ is a probability measure. In particular, the coincidence set $\bigcup_{j<k}\{\omega_j=\omega_k\}$, a countable union of hyperplanes, has measure zero for this density: $B$ has simple spectrum almost surely. Formula \eqref{EQ_weyl_integration} then expresses the law of $B=P\Delta P^\dagger$ as the pushforward of the product measure $dP\otimes(\tilde f\,d\omega)$ under the polar map $(P,\omega)\mapsto P\Delta P^\dagger$, i.e.
$$d\mu(B)=\tilde f(\omega)\,d\omega\otimes dP,$$
where $dP$ is the normalized Haar measure on the flag manifold $U(d)/\mathbb{T}^d$. Hence $P$ and $\omega$ are independent, with $P$ Haar-distributed.
\end{proof}

Geometrically, Proposition \ref{PRO_factorisation_mesure_Haar} says that the eigenbasis $P$ of $B$ is Haar-distributed (a random direction with no preferred orientation relative to that of the drift) and independent of $\omega$, whatever the isotropic ensemble. This decoupling is what makes the framework convenient: the bounds of Proposition \ref{PRO_rappel_estimations_tauSS} only involve $\tilde{A}^{(d)}=P^\dagger A^{(d)}P$, while $\omega$ enters only through the rational independence condition, which holds almost surely (as shown below); hence the statistical behavior of $\tau_{SS}^{(d)}$ does not depend on the particular ensemble from which $B$ is drawn.

\begin{proposition}[Almost Sure Rational Independence of the eigenvalues]\label{PRO_independence_rationnelle_generique}
Let the eigenvalue spectrum $\omega = (\omega_1, \dots, \omega_d) \in \mathbb{R}^d$ be a random vector whose probability distribution is absolutely continuous with respect to the Lebesgue measure, i.e., it admits a joint density $g \in L^1(\mathbb{R}^d)$. Then, the eigenvalues are rationally linearly independent with probability $1$.
\end{proposition}

\begin{proof}
The eigenvalues $\omega_1, \dots, \omega_d$ are rationally linearly dependent if and only if there exists a non-zero integer vector $n = (n_1, \dots, n_d) \in \mathbb{Z}^d \setminus \{0\}$ such that the inner product $n \cdot \omega = 0$.

For any fixed non-zero vector $n$, the set of solutions $H_n = \{ x \in \mathbb{R}^d \mid n \cdot x = 0 \}$ defines a hyperplane of dimension $d-1$ in $\mathbb{R}^d$. The $d$-dimensional Lebesgue measure of any such hyperplane is zero.

The set of all rationally dependent spectra, denoted $\mathcal{R}$, is the union of all these hyperplanes over the integer lattice:

$$\mathcal{R} = \bigcup_{n \in \mathbb{Z}^d \setminus \{0\}} H_n.$$

Since the set $\mathbb{Z}^d \setminus \{0\}$ is countably infinite, $\mathcal{R}$ is a countable union of sets of Lebesgue measure zero. By the countable subadditivity of measures, the Lebesgue measure of $\mathcal{R}$ is zero.

Because the probability distribution of $\omega$ admits a density $g$ with respect to the Lebesgue measure, the probability of drawing a rationally dependent spectrum is given by the integral of the density over this null set:

$$\mathbb{P}(\omega \in \mathcal{R}) = \int_{\mathcal{R}} g(\omega) d\omega = 0.$$

Consequently, the complementary event that the eigenvalues are rationally linearly independent occurs with probability $1$. 
\end{proof}
As a consequence of Proposition \ref{PRO_independence_rationnelle_generique}, the choice of $f$ in the factorization formula of Proposition \ref{PRO_factorisation_mesure_Haar} is irrelevant for the application of Proposition \ref{PRO_rappel_estimations_tauSS}: with probability 1, the rational rank of the eigenvalues of $B$ is maximal. (To avoid confusion with the spectrum $(\lambda_m)$ of the drift, the eigenvalues of the control operator $B$ are denoted by $\omega_1,\dots,\omega_d$.)

\begin{remark}\em
The law of the interaction matrix $\tilde{A}^{(d)} = P^\dagger A^{(d)} P$ only depends on the multiset $\{\lambda_1,\dots,\lambda_d\}$, and not on the order in which the eigenvalues of the drift are listed. Indeed, if $Q$ is a permutation matrix, $Q A^{(d)} Q^\dagger$ is the diagonal matrix with permuted entries, and
\[
P^\dagger \big(Q A^{(d)} Q^\dagger\big) P = (Q^\dagger P)^\dagger A^{(d)} (Q^\dagger P),
\]
where $Q^\dagger P$ is again Haar distributed. Note that such a permutation acts on the rows of $P$, that is, on the eigenbasis of the drift; the vertices of the interaction graph, which are indexed by the eigenvectors of $B$ (the columns of $P$), are not affected. Since the quantities $v_d$, $\bar{\lambda}_d$ and the Lyapunov ratio are also symmetric functions of the eigenvalues, the ordering $|\lambda_1| \le \dots \le |\lambda_d|$ used in Assumptions \ref{ASS_geometric_class}, \ref{ASS_spectral_dominance} and \ref{ASS_quantitative_dominance} is only a convention to identify the largest eigenvalues. When the moduli of a given sequence $(\lambda_m)_{m \geq 1}$ fail to be non-decreasing (as for $\lambda_m = m^a$ with $a<0$), the results apply to the family $\{\lambda_1,\dots,\lambda_d\}$ sorted, for each $d$, by increasing moduli; this reindexation may depend on $d$.
\end{remark}

\subsection{Spectrum-Independent Matrix Moments} \label{SEC_preliminary_moments}

The off-diagonal entries of the transformed matrix from the drift in the eigenbasis of $B$ are given by:
$$\tilde{A}_{j,k}^{(d)} = i \sum_{m=1}^d \lambda_m \overline{P}_{m,j} P_{m,k}.$$

For the sake of readability, we do the proof only when $j=1$ (the arguments are similar for other values of $j$).  This paragraph is  essentially a rephrasing  of Lemma 2.22 in \cite{Meckes_2019} using the Weingarten formalism (see Section 2 of \cite{Collins2003}).

Because $P$ is Haar-distributed, its columns are orthonormal vectors uniformly distributed on the unit sphere in $\mathbb{C}^d$. For $k \neq 1$, 
the matrix element after the change of basis is expressed as:
\begin{equation}\label{EQ_Expression_Ajk}
\tilde{A}_{1,k}^{(d)} = i \sum_{m=1}^d \lambda_m \overline{P}_{m,1} P_{m,k}
\end{equation}
We will now compute the expected value and the variance of $\tilde{A}_{1,k}^{(d)}$.
\begin{lemma}
For every $d$ in $\mathbb{N}$ and $1 < k\leq d$,
$ \displaystyle \mathbb{E}[\tilde{A}_{1,k}^{(d)}] =0 $.
\end{lemma}

\begin{proof}
The Haar measure is invariant under right-multiplication by any unitary matrix. Let us consider the 
diagonal unitary matrix $L$ defined by $L_{j,j} = 1$ for all $j \neq k$, and $L_{k,k} = -1$.
The transformation $P \mapsto P' = P L$ preserves the measure. The components of the modified matrix $P'$ satisfy $P'_{m,1} = P_{m,1}$ (since the column index $1$ is distinct from $k$) and $P'_{m,k} = -P_{m,k}$.
The expectation of the product must remain invariant under this gauge transformation:
$$\mathbb{E}[\overline{P}_{m,1} P_{m,k}] = \mathbb{E}[\overline{P'}_{m,1} P'_{m,k}]
= \mathbb{E}[\overline{P}_{m,1} (-P_{m,k})]  = - \mathbb{E}[\overline{P}_{m,1} P_{m,k}] .$$
This symmetry enforces $\mathbb{E}[\overline{P}_{m,1} P_{m,k}] = 0$ for every row index $m$. Since the integral of each term in the sum \eqref{EQ_Expression_Ajk} vanishes, the result
$$\mathbb{E}[\tilde{A}_{1,k}^{(d)}] = i \sum_{m=1}^d \lambda_m \mathbb{E}[\overline{P}_{m,1} P_{m,k}]= 0,$$ 
follows by linearity.
\end{proof}

\begin{lemma}\label{LEM_variance}
For every $d$ in $\mathbb{N}$ and $1 < k\leq d$,
 $ \displaystyle \mathbb{E}[|\tilde{A}_{1,k}^{(d)}|^2]= \dfrac{d}{d^2-1} v_d$.
\end{lemma}

To compute the variance of $\tilde{A}_{1,k}^{(d)}$, we
expand the modulus square of \eqref{EQ_Expression_Ajk} to get
\begin{equation}\label{EQ_module_carre_A_tilde}
|\tilde{A}_{1,k}^{(d)}|^2 = \sum_{m=1}^d \lambda_m^2 |P_{m,1}|^2 |P_{m,k}|^2 + \sum_{m \neq n} \lambda_m \lambda_n P_{m,k} P_{n,1} \overline{P}_{m,1} \overline{P}_{n,k}.
\end{equation}
We will compute successively the integrals of the square terms $|P_{m,1}|^2 |P_{m,k}|^2$ and of the cross-terms $P_{m,k} P_{n,1} \overline{P}_{m,1} \overline{P}_{n,k}$.

\begin{lemma}[Square terms Expectation]\label{LEM_Square_Terms}
$\mathbb{E}\left[|P_{m,1}|^2 |P_{m,k}|^2\right] =1/(d^2+d)$.
\end{lemma}
\begin{proof}
Let $P \in U(d)$ be a random unitary matrix distributed according to the Haar measure and let $k \neq 1$. We wish to compute the expectation value
$$\mathbb{E}\left[|P_{m,1}|^2 |P_{m,k}|^2\right] = \mathbb{E}\left[P_{m,1} P_{m,k} \overline{P}_{m,1} \overline{P}_{m,k}\right].$$
The general Weingarten formula for the integration of polynomials of matrix elements over the unitary group $U(d)$ reads
\begin{equation}
\mathbb{E}\left[P_{i_1, j_1} \dots P_{i_q, j_q} \overline{P}_{i'_1, j'_1} \dots \overline{P}_{i'_q, j'_q}\right] = \sum_{\sigma, \tau \in S_q} \delta_{\vec{i}, \sigma(\vec{i'})}\, \delta_{\vec{j}, \tau(\vec{j'})}\, \mathrm{Wg}(\sigma \tau^{-1}, d),
\end{equation}
where $S_q$ is the symmetric group of degree $q$, $\mathrm{Wg}$ is the Weingarten function (see \cite{Collins2003}), and $\delta_{\vec{i}, \sigma(\vec{i'})}$ is a shorthand for the product of Kronecker deltas $\prod_{a=1}^q \delta_{i_a, i'_{\sigma(a)}}$.

In our case, we evaluate a moment of degree $q=2$, and the symmetric group $S_2$ consists of the identity $e$ and the transposition $(12)$. The row and column indices are $i_1 = i_2 = m$, $j_1 = 1$, $j_2 = k$ for the unconjugated terms $P_{m,1}, P_{m,k}$, and $i'_1 = i'_2 = m$, $j'_1 = 1$, $j'_2 = k$ for the conjugated terms $\overline{P}_{m,1}, \overline{P}_{m,k}$. The integration formula therefore reduces to
\begin{equation}\label{EQ_Weingarten_q2_termes_carres}
\mathbb{E}\left[P_{m,1} P_{m,k} \overline{P}_{m,1} \overline{P}_{m,k}\right] = \sum_{\sigma, \tau \in S_2} \left( \delta_{m, m} \delta_{m, m} \right)_\sigma \left( \delta_{j_1, j'_{\tau(1)}} \delta_{j_2, j'_{\tau(2)}} \right)_\tau \mathrm{Wg}(\sigma \tau^{-1}, d).    
\end{equation}

We evaluate the constraints for both permutations in $S_2$. Since $i_1 = i_2 = i'_1 = i'_2 = m$, the row constraint is trivially satisfied: $\delta_{m, m} \delta_{m, m} = 1$ both for $\sigma = e$ and for $\sigma = (12)$. For the column constraint, with $j_1 = j'_1 = 1$ and $j_2 = j'_2 = k$, we get $\delta_{1, 1} \delta_{k, k} = 1$ for $\tau = e$, whereas $\delta_{1, k} \delta_{k, 1} = 0$ for $\tau = (12)$, because $k \neq 1$. The double sum over $S_2 \times S_2$ in \eqref{EQ_Weingarten_q2_termes_carres} thus collapses to the terms where $\tau = e$, and the sum over $\sigma \in \{e, (12)\}$ remains:
\[\mathbb{E}\left[|P_{m,1}|^2 |P_{m,k}|^2\right] = \sum_{\sigma \in \{e, (12)\}} \mathrm{Wg}(\sigma, d) = \mathrm{Wg}(e, d) + \mathrm{Wg}((12), d).\]

For $q=2$, the values of the Weingarten function are rational functions of the dimension $d$, namely $\mathrm{Wg}(e, d) = \frac{1}{d^2 - 1}$ and $\mathrm{Wg}((12), d) = -\frac{1}{d(d^2 - 1)}$ (see \cite{collins2006integration}, example at the end of Section 2.2). Substituting these values, we obtain
\[\mathbb{E}\left[|P_{m,1}|^2 |P_{m,k}|^2\right] = \frac{1}{d^2 - 1} - \frac{1}{d(d^2 - 1)} = \frac{1}{d(d+1)}.\]
\end{proof}

\begin{lemma}[Cross-term Expectation]\label{LEM_Cross_Term}
$\mathbb{E}\left[ P_{m,k} P_{n,1} \overline{P}_{m,1} \overline{P}_{n,k} \right] = 
-1/(d^3-d)$.
\end{lemma}
\begin{proof}
Let $P \in U(d)$ be a random unitary matrix distributed according to the normalized Haar measure. We evaluate the expectation value of the cross terms 
$$\mathbb{E}\left[ P_{m,k} P_{n,1} \overline{P}_{m,1} \overline{P}_{n,k} \right],$$
under the assumptions that $m \neq n$ and $k \neq 1$. As in the previous computation, we apply the Weingarten formula  (\cite{Weingarten1978, Collins2003})  in the permutation group $S_2$:
\begin{equation}\label{EQ_Weingarten_termes_croises}
\mathbb{E}\left[ P_{i_1, j_1} P_{i_2, j_2} \overline{P}_{i'_1, j'_1} \overline{P}_{i'_2, j'_2} \right] = \sum_{\sigma, \tau \in S_2} \delta_{\vec{i}, \sigma(\vec{i'})}\, \delta_{\vec{j}, \tau(\vec{j'})}\, \mathrm{Wg}(\sigma \tau^{-1}, d).
\end{equation}
The indices in \eqref{EQ_Weingarten_termes_croises} are $i_1 = m$, $i_2 = n$, $j_1 = k$, $j_2 = 1$ for the unconjugated elements $P_{m,k}, P_{n,1}$, and $i'_1 = m$, $i'_2 = n$, $j'_1 = 1$, $j'_2 = k$ for the conjugated elements $\overline{P}_{m,1}, \overline{P}_{n,k}$.

We evaluate the Kronecker delta constraints for $\sigma, \tau \in \{e, (12)\}$. The row constraint in
\eqref{EQ_Weingarten_termes_croises} 
$\delta_{i_1, i'_{\sigma(1)}} \delta_{i_2, i'_{\sigma(2)}}$ equals $\delta_{m, m} \delta_{n, n} = 1$ for $\sigma = e$, and $\delta_{m, n} \delta_{n, m} = 0$ for $\sigma = (12)$, since $m \neq n$; hence only the term $\sigma = e$ survives. The column constraint $\delta_{j_1, j'_{\tau(1)}} \delta_{j_2, j'_{\tau(2)}}$ equals $\delta_{k, 1} \delta_{1, k} = 0$ for $\tau = e$, since $k \neq 1$, and $\delta_{k, k} \delta_{1, 1} = 1$ for $\tau = (12)$; hence only the term $\tau = (12)$ survives.

The double sum in \eqref{EQ_Weingarten_termes_croises} therefore reduces to the single contribution of the pair $\sigma = e$, $\tau = (12)$. Since a transposition is its own inverse, $\mathrm{Wg}(e \cdot (12)^{-1}, d) = \mathrm{Wg}((12), d)$, and using the value $\mathrm{Wg}((12), d) = -\frac{1}{d(d^2 - 1)}$ of the Weingarten function for $q=2$, we obtain
\[\mathbb{E}\left[ P_{m,k} P_{n,1} \overline{P}_{m,1} \overline{P}_{n,k} \right] = -\frac{1}{d(d^2 - 1)},\]
which is the announced expression.
\end{proof}

\begin{proof}[Proof of Lemma \ref{LEM_variance}]
Injecting the respective expectations $\frac{1}{d(d+1)}$ (Lemma \ref{LEM_Square_Terms}, square terms) and $-\frac{1}{d(d^2-1)}$ (Lemma \ref{LEM_Cross_Term}, cross terms) in the expansion \eqref{EQ_module_carre_A_tilde} of $|\tilde{A}_{1,k}^{(d)}|^2$, 
the total expectation becomes:
$$\mathbb{E}[|\tilde{A}_{1,k}^{(d)}|^2] = \frac{1}{d(d+1)} \sum_{m=1}^d \lambda_m^2 - \frac{1}{d(d^2-1)} \sum_{m \neq n} \lambda_m \lambda_n.$$
The double sum over distinct indices can be rewritten in terms of simple sums:
$$\sum_{m \neq n} \lambda_m \lambda_n = \left( \sum_{m=1}^d \lambda_m \right)^2 - \sum_{m=1}^d \lambda_m^2.$$
Substituting this identity, we group the terms involving $\sum_m \lambda_m^2$:
$$\mathbb{E}[|\tilde{A}_{1,k}^{(d)}|^2] = \left( \frac{1}{d(d+1)} + \frac{1}{d(d^2-1)} \right) \sum_{m=1}^d \lambda_m^2 - \frac{1}{d(d^2-1)} \left( \sum_{m=1}^d \lambda_m \right)^2.$$
The prefactor simplifies as follows:
$$\frac{1}{d(d+1)} + \frac{1}{d(d+1)(d-1)} = \frac{d-1+1}{d(d+1)(d-1)} = \frac{1}{d^2-1}.$$

The expression of the variance is therefore:
$$\mathbb{E}[|\tilde{A}_{1,k}^{(d)}|^2] = \frac{1}{d^2-1} \sum_{m=1}^d \lambda_m^2 - \frac{1}{d(d^2-1)} \left( \sum_{m=1}^d \lambda_m \right)^2=
\frac{d}{d^2-1} v_d,
$$
which concludes the proof of Lemma \ref{LEM_variance}
\end{proof}

\begin{remark}\label{REM_variance_general_weights}\em
The computation above uses neither the ordering nor any growth property of the real numbers $\lambda_1,\dots,\lambda_d$. By invariance of the Haar measure under permutations of the columns, it holds for any pair $j \neq k$. Hence, for any real vector $\nu = (\nu_1,\dots,\nu_d)$ and any $j \neq k$,
\begin{equation}\label{EQ_variance_general_weights}
\mathbb{E}\Big[\Big|\sum_{m=1}^d \nu_m \overline{P}_{m,j} P_{m,k}\Big|^2\Big] = \frac{1}{d^2-1}\Big(\sum_{m=1}^d \nu_m^2 - \frac{1}{d}\Big(\sum_{m=1}^d \nu_m\Big)^2\Big) \leq \frac{1}{d^2-1}\sum_{m=1}^d \nu_m^2.
\end{equation}
Formula \eqref{EQ_variance_general_weights} also shows that the products $Y_m = d\, \overline{P}_{m,j} P_{m,k}$, $m=1,\dots,d$, are \emph{not} independent: they satisfy the constraint $\sum_m Y_m = d \langle P_{\cdot,j}, P_{\cdot,k}\rangle = 0$, and the cross-term computation gives $\mathbb{E}[Y_m \overline{Y_n}] = -\frac{d}{d^2-1}$ for $m\neq n$. These correlations are individually of order $1/d$, but their cumulative effect on $\sum_m \nu_m Y_m$ is of order $1$ unless $(\sum_m \nu_m)^2 = o(d \sum_m \nu_m^2)$. In particular, sums of the form $\sum_m \nu_m Y_m$ cannot be treated with central limit theorems for independent variables.
\end{remark}

\subsection{Spectrum-Independent Bound on MST}\label{SEC_Uniform_bound_MST}

 \begin{theorem}\label{TH_Uniform_bound_MST}
Let $A \in \mathbb{C}^{d \times d}$ be a diagonal drift matrix defined by $A = i \, 
\mathrm{diag}(\lambda_1, \dots, \lambda_d)$. Define the empirical spectral variance as 
$v_d = \frac{1}{d}\sum_{k=1}^d \lambda_k^2 - (\frac{1}{d}\sum_{k=1}^d \lambda_k)^2$, 
and assume $v_d > 0$. 

Let $P \in U(d)$ be a random unitary matrix distributed according to the Haar measure. Let $K_d$ denote the complete graph on $d$ vertices where each edge $(j,k)$ is assigned the weight $W_{jk} = |(P^\dagger A P)_{jk}|^{-1}$ for $1 \le j < k \le d$.

If $M_d$ represents the total weight of the Minimal Spanning Tree (MST) of $K_d$, then 
$M_d = \mathcal{O}_P(d^{3/2} v_d^{-1/2})$, or in other words, for every sequence $\kappa_d \to \infty$, $\mathbb{P}(M_d \ge \kappa_d d^{3/2} v_d^{-1/2}) \to 0$.
\end{theorem} 

\begin{proof}
Since the Minimal Spanning Tree (MST) is the spanning tree of minimal total weight, its weight $M_d$ is
bounded above by the weight of any arbitrary deterministic spanning tree. We select the path graph 
$T_{\mathrm{path}}$ defined by the edge set $\{(j, j+1) \mid 1 \le j \le d-1\}$. The weight satisfies:
\begin{equation}
M_d \le \sum_{j=1}^{d-1} |\tilde{A}_{j,j+1}|^{-1}
\end{equation}
where $\tilde{A} = P^\dagger A P$. By the invariance of the Haar measure under right multiplication by permutation matrices (which permutes simultaneously the rows and the columns of $\tilde{A}$), all off-diagonal elements of $\tilde{A}$ are identically distributed. Thus, the expected upper bound reduces to $\mathbb{E}[M_d] \le (d-1) \mathbb{E}[\vert{}\tilde{A}_{12}\vert{}^{-1}]$.

Let $U$ and $V$ denote the first two columns of the unitary matrix $P$. Conditionally on $U$, the vector $V$ is uniformly distributed on the unit sphere of the $(d-1)$-dimensional subspace $U^\perp$. We employ the generative representation $V \stackrel{d}{=} \Pi_{U^\perp} G / \Vert{}\Pi_{U^\perp} G\Vert{}$, where $G \sim \mathcal{CN}(0, I_d)$ is a standard isotropic complex Gaussian vector and $\Pi_{U^\perp} = I - UU^\dagger$ is the orthogonal projector onto $U^\perp$ (see for instance Chapter 2 of \cite{Meckes_2019}).

The target matrix element is the bilinear form $\tilde{A}_{12} = U^\dagger A V$. We define the projected vector $\mathbf{w} = \Pi_{U^\perp} A U = A U - (U^\dagger A U)U$. Its squared norm represents the spectral variance of $A$ along the state $U$, explicitly given by $V_A(U) = \Vert{}AU\Vert{}^2 - \vert{}U^\dagger A U\vert{}^2$ (Pythagorean theorem).
The magnitude of the matrix element can then be expressed as:
\begin{equation}
|\tilde{A}_{12}| = \sqrt{V_A(U)} \frac{|Z|}{\sqrt{|Z|^2 + R^2}}
\end{equation}
where $Z = \langle \mathbf{w}/\Vert{}\mathbf{w}\Vert{}, \Pi_{U^\perp} G \rangle \sim \mathcal{CN}(0,1)$ and $R := \sqrt{\Vert{}\Pi_{U^\perp} G\Vert{}^2 - \vert{}Z\vert{}^2} \ge 0$ (so that $R^2 \sim \Gamma(d-2, 1)$; for $d = 2$, $R \equiv 0$ and the estimate below is a fortiori valid) are independent random variables.

Using the triangle inequality $\sqrt{\vert{}Z\vert{}^2 + R^2} \le \vert{}Z\vert{} + R$, we derive an upper bound for the inverse magnitude:
\begin{equation}
|\tilde{A}_{12}|^{-1} \le \frac{1}{\sqrt{V_A(U)}} \left( 1 + \frac{R}{|Z|} \right)
\end{equation}
Taking the conditional expectation with respect to $U$, and exploiting the independent moments $\mathbb{E}[\vert{}Z\vert{}^{-1}] = \sqrt{\pi}$ and $\mathbb{E}[R] \le \sqrt{d-2}$, we obtain:
\begin{equation}
\mathbb{E}[|\tilde{A}_{12}|^{-1} \mid U] \le \frac{1 + \sqrt{\pi(d-2)}}{\sqrt{V_A(U)}}
\end{equation}

To compute the unconditional expectation, we must evaluate $\mathbb{E}[V_A(U)^{-1/2}]$. The algebraic expansion of the variance is $V_A(U) = \sum_{1 \le i < j \le d} (\lambda_i - \lambda_j)^2 \vert{}U_i\vert{}^2 \vert{}U_j\vert{}^2$.
Because $U$ is uniformly distributed on the complex unit sphere, the vector of its squared moduli $(\vert{}U_1\vert{}^2, \dots, \vert{}U_d\vert{}^2)$ follows a uniform Dirichlet distribution on the standard $(d-1)$-simplex. This vector admits the stochastic representation $\vert{}U_i\vert{}^2 \stackrel{d}{=} X_i/Y$, where $X_1, \dots, X_d$ are independent standard exponential variables ($X_i \sim \mathrm{Exp}(1)$) and $Y = \sum_{k=1}^d X_k$ acts as the global scaling variable following a Gamma distribution, $Y \sim \Gamma(d,1)$ (see for instance Chapter 2 of \cite{Meckes_2019}).

The sum $Y$ and the normalized vector $(X_1/Y, \dots,$  $X_d/Y)$ on the simplex are statistically independent. We define the unnormalized quadratic form $D(X) = \sum_{1 \le i < j \le d} (\lambda_i - \lambda_j)^2 X_i X_j$.
By substituting the generative representation, the conditional variance scales as:
\begin{equation}
V_A(U) \stackrel{d}{=} \frac{1}{Y^2} D(X)
\end{equation}
Taking the inverse square root yields $D(X)^{-1/2} \stackrel{d}{=} Y^{-1} V_A(U)^{-1/2}$.
Since $V_A(U)$ is a function of the squared moduli $(\vert{}U_1\vert{}^2, \dots, \vert{}U_d\vert{}^2) = (X_1/Y, \dots, X_d/Y)$, which are independent of $Y$ (as noted above), the expectation factorizes:
\begin{equation}
\mathbb{E}[D(X)^{-1/2}] = \mathbb{E}[Y^{-1}] \mathbb{E}[V_A(U)^{-1/2}].
\end{equation}
The inverse moment of the Gamma distribution evaluates as:
\begin{equation}
\mathbb{E}[Y^{-1}] = \int_0^\infty y^{-1} \frac{y^{d-1} e^{-y}}{\Gamma(d)} dy = \frac{\Gamma(d-1)}{\Gamma(d)} = \frac{1}{d-1}
\end{equation}
Rearranging the terms isolates the target expectation:
\begin{equation}
\mathbb{E}[V_A(U)^{-1/2}] = (d-1) \mathbb{E}\left[ \left( \sum_{1 \le i < j \le d} (\lambda_i - \lambda_j)^2 X_i X_j \right)^{-1/2} \right]
\end{equation}

We now apply Jensen's inequality to the strictly convex function $x \mapsto x^{-1/2}$. We define the nonnegative weights $w_{ij} = (\lambda_i - \lambda_j)^2$ and normalize them by their total sum $W = \sum_{i < j} w_{ij} = d^2 v_d$:
\begin{eqnarray}
\mathbb{E}[V_A(U)^{-1/2}] &\!=\!&\!\!\!\!(d-1) \mathbb{E} \left[ \left( W \!\!\!\!\sum_{1 \le i < j \le d} \frac{w_{ij}}{W} X_i X_j \right)^{-1/2} \right] =
\frac{d-1}{\sqrt{W}} 
\mathbb{E} \left[ \left( \sum_{1 \le i < j \le d} \frac{w_{ij}}{W} X_i X_j \right)^{-1/2} \right]
\nonumber \\
&\leq &
 \frac{d-1}{\sqrt{W}} \sum_{i < j} \frac{w_{ij}}{W} \mathbb{E}[(X_i X_j)^{-1/2}]
\end{eqnarray}
Using the independence of the exponential variables, the moment evaluates to $\mathbb{E}[(X_i X_j)^{-1/2}] = (\mathbb{E}[X_1^{-1/2}])^2 = (\Gamma(1/2))^2 = \pi$. Substituting $W$, we bound the expectation:
\begin{equation}
\mathbb{E}[V_A(U)^{-1/2}] \le \frac{(d-1)\pi}{d\sqrt{v_d}} < \frac{\pi}{\sqrt{v_d}}.
\end{equation}

Combining the estimates, the unconditional expectation of a generic inverse edge weight satisfies:
\begin{equation}
\mathbb{E}[|\tilde{A}_{12}|^{-1}] \le \frac{\pi + \pi^{3/2}\sqrt{d-2}}{\sqrt{v_d}} = \mathcal{O}\left( \sqrt{\frac{d}{v_d}} \right).
\end{equation}
Consequently, the expected total weight of the path graph is bounded 
$$\mathbb{E}[M_d] \! \le \! (d-1) \mathbb{E}[\vert{}\tilde{A}_{12}\vert{}^{-1}] = \mathcal{O}(d^{3/2} v_d^{-1/2}).$$
By Markov's inequality, for any divergent sequence $\kappa_d \to \infty$ as $d \to \infty$:
\begin{equation}
\mathbb{P}\left( M_d \ge \kappa_d \frac{d^{3/2}}{\sqrt{v_d}} \right) \le \frac{\mathbb{E}[M_d]}{\kappa_d d^{3/2} v_d^{-1/2}} \to 0.
\end{equation}
This establishes the stochastic domination $M_d = \mathcal{O}_P(d^{3/2} v_d^{-1/2})$, concluding the proof.
\end{proof}

\begin{remark}\label{REM_explicit_upper_bound}\em
The proof gives the explicit, non-asymptotic estimate $\mathbb{E}[M_d] \le (d-1)\frac{\pi + \pi^{3/2}\sqrt{d-2}}{\sqrt{v_d}} \le (\pi + \pi^{3/2}) \frac{d^{3/2}}{\sqrt{v_d}}$. Combined with Proposition \ref{PRO_rappel_estimations_tauSS} and Markov's inequality, it yields, for every spectrum with $v_d>0$, every $d$ and every $\epsilon \in (0,1)$:
\[
\mathbb{P}\left( \tau_{SS}^{(d)} \le C_\epsilon \frac{d^{3/2}}{\sqrt{v_d}} \right) \geq 1 - \epsilon, \qquad C_\epsilon = \frac{\pi(\pi + \pi^{3/2})}{2\epsilon}.
\]
\end{remark}

We now give the counterpart of Theorem \ref{TH_Uniform_bound_MST} for the lower bound. It is crude (compare with Theorem \ref{TH_coupe_minimale_Lyapunov} and Proposition \ref{PRO_lower_bound_dominance}), but it holds for any spectrum.

\begin{proposition}[Spectrum-independent lower bound]\label{PRO_Uniform_lower_bound}
Let $A = i \, \mathrm{diag}(\lambda_1, \dots, \lambda_d)$ with $v_d > 0$, let $P \in U(d)$ be Haar distributed and $C_1 = \sum_{k \neq 1} |\tilde{A}_{1k}|^2$ be the capacity of the star cut isolating the vertex $1$. Then $\mathbb{E}[C_1] = \frac{d}{d+1} v_d$ and, for every $\epsilon \in (0,1)$ and every $d$,
\[
\mathbb{P}\left(\tau_{SS}^{(d)} \geq \frac{2\sqrt{\epsilon}}{\sqrt{v_d}}\right) \geq 1 - \epsilon.
\]
\end{proposition}

\begin{proof}
By Lemma \ref{LEM_variance}, $\mathbb{E}[C_1] = (d-1)\frac{d}{d^2-1} v_d = \frac{d}{d+1} v_d \le v_d$. By Markov's inequality, $\mathbb{P}(C_1 \geq v_d/\epsilon) \leq \epsilon$. On the complementary event, Proposition \ref{PRO_rappel_estimations_tauSS} gives $\tau_{SS}^{(d)} \geq 2/\sqrt{C_1} \geq 2\sqrt{\epsilon}/\sqrt{v_d}$.
\end{proof}

\subsection{Gaussian approximation of the entries of Haar matrices}\label{SEC_Gaussian_Approximation}

\begin{theorem}[Borel-Poincar\'e lemma for complex spheres]\label{THE_Borel_Poincare}
Let $U^{(d)}$ be a random vector uniformly distributed on the unit sphere $S_{\mathbb{C}}^{d-1}$ of $\mathbb{C}^d$, and let $k \geq 1$ be a fixed integer. Then, as $d \to \infty$,
$$\sqrt{d}\, \big( U_1^{(d)}, \dots, U_k^{(d)} \big) \xrightarrow{\mathcal{L}} (Z_1, \dots, Z_k),$$
where $Z_1, \dots, Z_k$ are independent $\mathcal{CN}(0,1)$ variables.
\end{theorem}

\begin{proof}
Let $G \sim \mathcal{CN}(0, I_d)$. Since the law of $G$ is invariant under $U(d)$ and $G \neq 0$ almost surely, the vector $G/\|G\|$ is uniformly distributed on $S_{\mathbb{C}}^{d-1}$. Hence
$$\sqrt{d}\, \big( U_1^{(d)}, \dots, U_k^{(d)} \big) \stackrel{\mathcal{L}}{=} \frac{1}{\sqrt{S_d}}\, (G_1, \dots, G_k), \qquad S_d = \frac{1}{d} \sum_{m=1}^d |G_m|^2.$$
The variables $|G_m|^2$ are independent with mean $1$ and variance $1$, so that $S_d \xrightarrow{\mathbb{P}} 1$ by the weak law of large numbers, and $1/\sqrt{S_d} \xrightarrow{\mathbb{P}} 1$ by the continuous mapping theorem. Since the law of $(G_1, \dots, G_k)$ does not depend on $d$, Slutsky's theorem concludes the proof.
\end{proof}

Since a column of a Haar distributed matrix $P$, and the conjugate of a row of $P$ (which is a column of the Haar distributed matrix $P^\dagger$), are uniformly distributed on $S_{\mathbb{C}}^{d-1}$, and since the Haar measure is invariant under permutations of rows and columns, Theorem \ref{THE_Borel_Poincare} applies to any fixed number of entries of a given row or column of $P$. The following lemma extends it to entries taken in two columns.

\begin{lemma}\label{LEM_block_Gaussian}
Let $P$ be a Haar distributed matrix in $U(d)$, let $p \geq 1$ be a fixed integer, let $i_1,\dots,i_p$ be distinct row indices and $j \neq k$ two column indices (all possibly depending on $d$). Then, as $d \to \infty$, the random vector
$$\sqrt{d}\,\big(P_{i_1 j}, \dots, P_{i_p j}, P_{i_1 k}, \dots, P_{i_p k}\big),$$
converges in distribution to a vector of $2p$ independent $\mathcal{CN}(0,1)$ variables.
\end{lemma}

\begin{proof}
By invariance of the Haar measure under permutations of rows and columns, we may assume that $(i_1,\dots,i_p) = (1,\dots,p)$ and $(j,k) = (1,2)$. Let $g_1, g_2$ be independent $\mathcal{CN}(0,I_d)$ vectors, and consider their Gram-Schmidt orthonormalization
$$u_1 = \frac{g_1}{\|g_1\|}, \qquad h = g_2 - \frac{\langle g_1, g_2\rangle}{\|g_1\|^2}\, g_1, \qquad u_2 = \frac{h}{\|h\|}.$$
The law of $(g_1, g_2)$ is invariant under $U(d)$ and the Gram-Schmidt procedure commutes with the action of $U(d)$, so the law of $(u_1, u_2)$ is a $U(d)$-invariant probability measure on the set of orthonormal pairs of $\mathbb{C}^d$. Since $U(d)$ acts transitively on this set, this measure is unique, and $(u_1,u_2)$ has the same law as the first two columns of $P$.

By the law of large numbers, $\|g_1\|^2/d \xrightarrow{\mathbb{P}} 1$ and $\|g_2\|^2/d \xrightarrow{\mathbb{P}} 1$. Conditionally on $g_1$, $\langle g_1, g_2\rangle \sim \mathcal{CN}(0, \|g_1\|^2)$, so that $\langle g_1, g_2\rangle/\|g_1\|^2 = \mathcal{O}_{\mathbb{P}}(d^{-1/2})$ and $\|h\|^2 = \|g_2\|^2 - |\langle g_1, g_2\rangle|^2/\|g_1\|^2 = d\,(1+o_{\mathbb{P}}(1))$. Hence, for $1 \leq i \leq p$,
$$\sqrt{d}\, u_{1,i} = \frac{\sqrt{d}}{\|g_1\|}\, g_{1,i}, \qquad \sqrt{d}\, u_{2,i} = \frac{\sqrt{d}}{\|h\|}\Big( g_{2,i} - \frac{\langle g_1, g_2\rangle}{\|g_1\|^2}\, g_{1,i} \Big),$$
where the prefactors converge in probability to $1$ and the correction term $\frac{\langle g_1, g_2\rangle}{\|g_1\|^2}\, g_{1,i}$ converges in probability to $0$. Since the law of $(g_{1,1},\dots,g_{1,p},g_{2,1},\dots,g_{2,p})$ is that of $2p$ independent $\mathcal{CN}(0,1)$ variables and does not depend on $d$, Slutsky's theorem concludes the proof.
\end{proof}

\section{Asymptotic normality of the drift entries under the Lyapunov condition}\label{SEC_Lyapunov_normality}

In this section, we establish that, if the spectrum satisfies the Lyapunov condition (Assumption \ref{ASS_Lyapunov_condition}), then
the scaled off-diagonal elements of $\tilde{A}^{(d)}$ are asymptotically standard complex normal variables.

\begin{proposition}\label{PRO_off_diagonal_normality}
Assume that the spectrum of $A^{(d)}$ satisfies the Lyapunov condition (Assumption \ref{ASS_Lyapunov_condition}). 
Define the scaled off-diagonal elements of the interaction matrix as $Z_{j,k}^{(d)} = \sqrt{\frac{d}{v_d}} \tilde{A}_{j,k}^{(d)}$, $j \neq k$.
Then the law of $Z_{j,k}^{(d)}$ does not depend on the pair $(j,k)$ and, for any sequence of pairs $(j_d,k_d)$ with $j_d \neq k_d$, 
$$Z_{j_d,k_d}^{(d)} \xrightarrow{\mathcal{L}} \mathcal{CN}(0,1) \quad \mbox{ as } d \to \infty.$$
\end{proposition}

Proposition \ref{PRO_off_diagonal_normality} adapts a classical result by Meckes (Theorem 2.23 in \cite{Meckes_2019}), originally established for the linear functional $P \mapsto \mathrm{Tr}(A P)$, to the bilinear form $P \mapsto P^\dagger_{\cdot,j} A P_{\cdot,k}$. Following the probabilistic framework of \cite{Meckes_2019} (Chapter 2), our proof proceeds in several steps:

\begin{enumerate}
\item Using the Weingarten calculus, we computed the finite-dimensional first and second moments of the variables $Z_{j,k}^{(d)}$: $\mathbb{E}[Z_{j,k}^{(d)}] = 0$ and $\mathbb{E}[|Z_{j,k}^{(d)}|^2] = \frac{d^2}{d^2-1} \to 1$ (Section \ref{SEC_preliminary_moments}).
\item Conditioning on an arbitrary column of $P$ (e.g., $P_{\cdot,1}$), any other column is uniformly distributed on the unit sphere of the $(d-1)$-dimensional orthogonal complement. Using the standard generative representation by Gaussian vectors, we compute the conditional variance of the entries of $\tilde{A}$ with respect to $P_{\cdot,1}$. Under the Lyapunov condition, this conditional variance is asymptotically equivalent, in probability, to $v_d/d$. This is done in Section \ref{SEC_Conditional_variance}.
\item Combining this conditional structure with Slutsky's theorem, we prove the asymptotic normality of $Z_{j,k}^{(d)}$ (Section \ref{SEC_normality}). The argument keeps track of the variance and of the orthogonality of the column vectors, which is essential: the summands of $\tilde{A}_{j,k}^{(d)} = i\sum_m \lambda_m \overline{P}_{m,j}P_{m,k}$ are \emph{not} independent (Remark \ref{REM_variance_general_weights}).
\end{enumerate}

By invariance of the Haar measure under right multiplication by permutation matrices, the law of $Z_{j,k}^{(d)}$ does not depend on $(j,k)$, and we may restrict ourselves to $(j,k) = (1,k)$ with $k \neq 1$.

\subsection{Conditional variance}\label{SEC_Conditional_variance}

We already established that the off-diagonal element takes the form of the bilinear sum:
$$\tilde{A}_{1,k}^{(d)} = i\sum_{m=1}^d \lambda_m \overline{P}_{m,1} P_{m,k}.$$
To evaluate the asymptotic behavior of this sum, we proceed via a conditioning argument. Let us fix the first column $P_{\cdot, 1}$. Conditionally on this vector, the column $P_{\cdot, k}$ (for $k \neq 1$) is uniformly distributed on the unit sphere of the $(d-1)$-dimensional subspace orthogonal to $P_{\cdot, 1}$, denoted $P_{\cdot, 1}^\perp$.

With the convention $\langle x, y \rangle = \sum_m \overline{x}_m y_m$, the sum can be written as the inner product 
$\tilde{A}_{1,k}^{(d)} = \langle \bar{c}, P_{\cdot, k} \rangle$, where the vector $\bar{c}$ has components $\bar{c}_m = -i \lambda_m P_{m,1}$. Because $P_{\cdot, k}$ belongs to $P_{\cdot, 1}^\perp$, $\langle \bar{c}, P_{\cdot, k} \rangle = \langle \Pi \bar{c}, P_{\cdot, k} \rangle$, where $\Pi$ denotes the orthogonal projection onto $P_{\cdot, 1}^\perp$.

By the rotational invariance of the uniform measure on the unit sphere of $P_{\cdot, 1}^\perp$, the conditional variance of the inner product depends exclusively on the norm of the projected vector:
$$\text{Var}(\tilde{A}_{1,k}^{(d)} \mid P_{\cdot, 1}) = \frac{\|\Pi \bar{c}\|^2}{d-1} = \frac{\|\bar{c}\|^2 - |\langle P_{\cdot, 1}, \bar{c} \rangle|^2}{d-1}.$$
Since $\|\bar{c}\|^2 = \sum_{m} \lambda_m^2 |P_{m,1}|^2$ and $\langle P_{\cdot, 1}, \bar{c} \rangle = -i \sum_{m} \lambda_m |P_{m,1}|^2$, we obtain
\begin{equation}\label{EQ_conditional_variance}
\|\Pi \bar{c}\|^2 = \mathcal{V}_d(P_{\cdot, 1}) := \sum_{m=1}^d \lambda_m^2 |P_{m,1}|^2 - \left( \sum_{m=1}^d \lambda_m |P_{m,1}|^2 \right)^2, \qquad \text{Var}(\tilde{A}_{1,k}^{(d)} \mid P_{\cdot, 1}) = \frac{\mathcal{V}_d(P_{\cdot, 1})}{d-1}.
\end{equation}
The quantity $\mathcal{V}_d(P_{\cdot, 1})$ is the variance of the spectrum $(\lambda_m)$ with respect to the probability weights $(|P_{m,1}|^2)_m$. It is also the capacity $C_1$ of the star cut isolating the vertex $1$ (see the proof of Theorem \ref{TH_coupe_minimale_Lyapunov}).

\begin{lemma}[Concentration of the conditional variance under the Lyapunov condition]
\label{LEM_concentration}
Let $A^{(d)} = i \, \mathrm{diag}(\lambda_1, \dots, \lambda_d)$ be a sequence of deterministic drift matrices with empirical spectral variance $v_d > 0$, satisfying the Lyapunov condition (Assumption \ref{ASS_Lyapunov_condition}).
Let $P_{\cdot, 1}$ be a random vector uniformly distributed on the complex unit sphere $S_{\mathbb{C}}^{d-1}$. Then, as $d \to \infty$,
$$\frac{\mathcal{V}_d(P_{\cdot, 1})}{v_d} \xrightarrow{\mathbb{P}} 1.$$
\end{lemma}

\begin{proof}
Since $\sum_m |P_{m,1}|^2 = 1$, $\mathcal{V}_d(P_{\cdot, 1})$ is invariant under any translation $\lambda_m \mapsto \lambda_m - c$ of the spectrum. The same holds for $v_d$ and for Assumption \ref{ASS_Lyapunov_condition}. We may therefore assume without loss of generality that the spectrum is centered: $\bar{\lambda}_d = 0$, so that $\frac{1}{d}\sum_m \lambda_m^2 = v_d$.

The unit vector $P_{\cdot, 1}$ admits the stochastic representation $P_{m,1} \stackrel{\mathcal{L}}{=} X_m / \sqrt{\sum_{j=1}^d |X_j|^2}$, where the $X_m \sim \mathcal{CN}(0,1)$ are independent. The squared moduli $I_m = |X_m|^2$ are independent exponential random variables with parameter $1$ ($\mathbb{E}[I_m] = \text{Var}(I_m) = 1$). Define
$$S_d = \frac{1}{d}\sum_{m=1}^d I_m, \qquad N_2 = \frac{1}{d v_d}\sum_{m=1}^d \lambda_m^2 I_m, \qquad N_1 = \frac{1}{d \sqrt{v_d}}\sum_{m=1}^d \lambda_m I_m.$$
Substituting the representation, we get the identity in law
$$\frac{\mathcal{V}_d(P_{\cdot, 1})}{v_d} \stackrel{\mathcal{L}}{=} \frac{N_2}{S_d} - \frac{N_1^2}{S_d^2}.$$
By independence of the $I_m$:
\begin{itemize}
\item[($i$)] $\mathbb{E}[S_d] = 1$ and $\text{Var}(S_d) = 1/d$;
\item[($ii$)] $\mathbb{E}[N_2] = 1$ and $\text{Var}(N_2) = \dfrac{\frac{1}{d}\sum_m \lambda_m^4}{d v_d^2}$, which tends to $0$ by Assumption \ref{ASS_Lyapunov_condition} (recall that the spectrum is centered);
\item[($iii$)] $\mathbb{E}[N_1] = \frac{1}{d\sqrt{v_d}}\sum_m \lambda_m = 0$ and $\text{Var}(N_1) = \dfrac{\frac{1}{d}\sum_m \lambda_m^2}{d v_d} = \dfrac{1}{d}$.
\end{itemize}
By Chebyshev's inequality, $(S_d, N_2, N_1) \xrightarrow{\mathbb{P}} (1,1,0)$. Since the function $(s,n_2,n_1) \mapsto n_2/s - n_1^2/s^2$ is continuous at $(1,1,0)$, the continuous mapping theorem yields $\mathcal{V}_d(P_{\cdot, 1})/v_d \xrightarrow{\mathbb{P}} 1$.
\end{proof}

\subsection{Asymptotic normality of the distribution}\label{SEC_normality}

\begin{proof}[Proof of Proposition \ref{PRO_off_diagonal_normality}]
We consider $Z_{1,k}^{(d)} = \sqrt{\frac{d}{v_d}} \langle \bar{c}, P_{\cdot, k} \rangle$, $k \neq 1$, with the notations of Section \ref{SEC_Conditional_variance}.

Let $G \sim \mathcal{CN}(0, I_d)$ be a standard complex Gaussian vector in $\mathbb{C}^d$, independent of $P_{\cdot, 1}$. Conditionally on $P_{\cdot, 1}$, the vector $\Pi G$ is a standard complex Gaussian vector of the $(d-1)$-dimensional space $P_{\cdot, 1}^\perp$, and $\Pi G/\|\Pi G\|$ is uniformly distributed on its unit sphere. Hence the pair $(P_{\cdot,1}, P_{\cdot,k})$ has the same law as $(P_{\cdot,1}, \Pi G/\|\Pi G\|)$, and, $\Pi$ being self-adjoint and idempotent,
$$Z_{1,k}^{(d)} \stackrel{\mathcal{L}}{=} \sqrt{\frac{d}{v_d}} \frac{\langle \Pi \bar{c}, G \rangle}{\|\Pi G\|}.$$

Since $v_d>0$, the spectrum is not constant, and since all the coordinates of $P_{\cdot,1}$ are almost surely non-zero, $\|\Pi\bar{c}\|^2 = \mathcal{V}_d(P_{\cdot,1}) > 0$ almost surely. Define
$$\xi = \frac{\langle \Pi \bar{c}, G \rangle}{\|\Pi \bar{c}\|}.$$
Conditionally on $P_{\cdot, 1}$, $\Pi\bar{c}/\|\Pi \bar{c}\|$ is a fixed unit vector and $\xi \sim \mathcal{CN}(0,1)$. Since this conditional law does not depend on $P_{\cdot, 1}$, the variable $\xi$ is independent of $P_{\cdot, 1}$ and $\xi \sim \mathcal{CN}(0,1)$. Using \eqref{EQ_conditional_variance}, we obtain the factorization
$$Z_{1,k}^{(d)} \stackrel{\mathcal{L}}{=} \sqrt{\frac{\mathcal{V}_d(P_{\cdot,1})}{v_d}} \cdot \frac{\sqrt{d}}{\|\Pi G\|} \cdot \xi.$$
We analyze the three factors:
\begin{itemize}
\item[($i$)] By Lemma \ref{LEM_concentration}, $\sqrt{\mathcal{V}_d(P_{\cdot,1})/v_d} \xrightarrow{\mathbb{P}} 1$.
\item[($ii$)] Conditionally on $P_{\cdot,1}$, $\|\Pi G\|^2$ is the sum of the squared moduli of the $d-1$ coordinates of $G$ in an orthonormal basis of $P_{\cdot, 1}^\perp$, hence $\|\Pi G\|^2 \sim \Gamma(d-1, 1)$ (equivalently, $2\|\Pi G\|^2 \sim \chi^2_{2(d-1)}$). Its mean is $d-1$ and its variance is $d-1$, so that, by Chebyshev's inequality, $\|\Pi G\|^2/d \xrightarrow{\mathbb{P}} 1$ and $\sqrt{d}/\|\Pi G\| \xrightarrow{\mathbb{P}} 1$.
\item[($iii$)] $\xi \sim \mathcal{CN}(0,1)$ for every $d$.
\end{itemize}
By Slutsky's theorem the product of a sequence with fixed law $\mathcal{CN}(0,1)$ by sequences converging in probability to the constant $1$ converges in distribution to $\mathcal{CN}(0,1)$. Hence $Z_{1,k}^{(d)} \xrightarrow{\mathcal{L}} \mathcal{CN}(0,1)$, which concludes the proof, since the law of $Z_{j,k}^{(d)}$ does not depend on $(j,k)$.
\end{proof}

\begin{remark}\em
The joint asymptotic behavior of several off-diagonal entries (asymptotic independence of a finite family of entries) is not needed in the sequel and is not addressed here: the bounds on the control time only involve the marginal law of one entry (upper bound) and the law of one column of $P$ (lower bound).
\end{remark}

 \section{Time estimates}\label{SEC_Time_estimates}

\subsection{Time estimates under the Lyapunov condition}

In this paragraph, we prove Theorem \ref{THE_random_encadrement}. The right hand side is a consequence of Theorem \ref{TH_Uniform_bound_MST} (see Remark \ref{REM_explicit_upper_bound}). For the left hand side, we use Theorem \ref{TH_coupe_minimale_Lyapunov} below.

\begin{theorem}\label{TH_coupe_minimale_Lyapunov}
Let $A^{(d)} = i \, \mathrm{diag}(\lambda_1, \dots, \lambda_d)$ be a sequence of deterministic drift matrices with $v_d > 0$, satisfying the Lyapunov condition (Assumption \ref{ASS_Lyapunov_condition}).
Let $P \in U(d)$ be a random unitary matrix distributed according to the Haar measure, $D$ be a random diagonal matrix with purely imaginary entries distributed according to any law absolutely continuous with respect to the Lebesgue measure, and $B=P D P^\dagger $. 

Then, for every $\delta > 0$, the state-to-state control time of the system $\psi'=(A^{(d)}+u(t)B) \psi$ satisfies
$$\lim_{d \to \infty} \mathbb{P}\left( \tau_{SS}^{(d)} \geq \frac{2-\delta}{\sqrt{v_d}} \right) = 1.$$
\end{theorem}

\begin{proof}
By Proposition \ref{PRO_rappel_estimations_tauSS} (applied with $j=1$),
$$\tau_{SS}^{(d)} \ge \frac{2}{\sqrt{C_1}}, \qquad C_1 = \sum_{k=2}^d |\tilde{A}_{1k}|^2,$$
with $\tilde{A}= P^\dagger A P $. The quantity $C_1$ is the capacity of the star cut isolating the vertex $1$ from the rest of the interaction graph.

To evaluate this sum, we complete the row by adding and subtracting the squared modulus of the diagonal element:
$$C_1 = \sum_{k=1}^d |\tilde{A}_{1k}|^2 - |\tilde{A}_{11}|^2.$$
The full sum of the squared moduli of a row is the diagonal element of the product $\tilde{A} \tilde{A}^\dagger$. Writing $A = i\Lambda$ with $\Lambda = \mathrm{diag}(\lambda_1,\dots,\lambda_d)$, we have $A A^\dagger = \Lambda^2$. Therefore:
$$\sum_{k=1}^d |\tilde{A}_{1k}|^2 = (\tilde{A} \tilde{A}^\dagger)_{11} = (P^\dagger \Lambda^2 P)_{11} = \sum_{m=1}^d \lambda_m^2 |P_{m1}|^2.$$
The diagonal element $\tilde{A}_{11}$ is the mean of the drift spectrum weighted by the same coefficients:
$$\tilde{A}_{11} = (P^\dagger A P)_{11} = i \sum_{m=1}^d \lambda_m |P_{m1}|^2.$$
Hence $C_1 = \mathcal{V}_d(P_{\cdot,1})$, with the notation \eqref{EQ_conditional_variance}.

The vector $P_{\cdot, 1}$ is uniformly distributed on the complex unit sphere $S_{\mathbb{C}}^{d-1}$. By Lemma \ref{LEM_concentration}, $C_1/v_d \xrightarrow{\mathbb{P}} 1$. Let $\delta \in (0,2)$ and $\eta>0$ be such that $(1+\eta)^{-1/2} \geq 1-\delta/2$. Then
$$\mathbb{P}\left( \tau_{SS}^{(d)} < \frac{2-\delta}{\sqrt{v_d}} \right) \leq \mathbb{P}\left( \frac{2}{\sqrt{C_1}} < \frac{2}{\sqrt{(1+\eta) v_d}} \right) = \mathbb{P}\left( C_1 > (1+\eta) v_d \right) \xrightarrow[d\to\infty]{} 0,$$
which concludes the proof.
\end{proof}

\begin{proof}[Proof of Theorem \ref{THE_random_encadrement}]
Let $\epsilon \in (0,1)$. By Remark \ref{REM_explicit_upper_bound}, $\mathbb{P}(\tau_{SS}^{(d)} > C_{\epsilon/2} d^{3/2}/\sqrt{v_d}) \leq \epsilon/2$ for every $d$. By Theorem \ref{TH_coupe_minimale_Lyapunov} (with $\delta = \epsilon$), there exists $d_0$ such that $\mathbb{P}(\tau_{SS}^{(d)} < (2-\epsilon)/\sqrt{v_d}) \leq \epsilon/2$ for every $d \geq d_0$. The result follows from the union bound.
\end{proof}

\subsection{Time estimates for spectra in the geometric class}

In this section, we establish the proof of Theorem \ref{THE_random_encadrement_geometric} for spectra in the uniformly geometric class $\mathcal{G}(a,b)$ (Assumption \ref{ASS_geometric_class}). Both inequalities follow from the spectrum-independent estimates of Section \ref{SEC_Uniform_bound_MST}, once the asymptotic behavior of the empirical spectral variance $v_d$ is known (Lemma \ref{LEM_variance_geometric}).

\begin{lemma}\label{LEM_variance_geometric}
Assume the spectrum belongs to the uniformly geometric class $\mathcal{G}(a,b)$ (Assumption \ref{ASS_geometric_class}). Let $K_1 = 1 + \frac{1}{1-a^{-1}}$ and $K_2 = 1 + \frac{1}{1-a^{-2}}$. There exists $d_1 \in \mathbb{N}$ such that, for every $d \geq d_1$,
\[
\sum_{m=1}^d |\lambda_m| \leq K_1 |\lambda_d|, \qquad \sum_{m=1}^d \lambda_m^2 \leq K_2 \lambda_d^2, \qquad 
\left(1 - \frac{K_1^2}{d}\right)\frac{\lambda_d^2}{d} \leq v_d \leq K_2 \frac{\lambda_d^2}{d}.
\]
In particular, $v_d \asymp \lambda_d^2/d$ as $d \to \infty$.
\end{lemma}

\begin{proof}
For $m_0 \leq d-k \leq d$, iterating the inequality $|\lambda_{m+1}| \geq a |\lambda_m|$ gives $|\lambda_{d-k}| \le a^{-k} |\lambda_d|$. For $m < m_0$, $|\lambda_m| \leq |\lambda_{m_0}| \leq a^{-(d-m_0)} |\lambda_d|$, so that $\sum_{m<m_0} |\lambda_m| \leq m_0 a^{-(d-m_0)} |\lambda_d| \le |\lambda_d|$ and $\sum_{m<m_0} \lambda_m^2 \leq \lambda_d^2$ for $d$ large enough. Summing the geometric series over the indices $m \geq m_0$ yields the first two inequalities. 

For the last one, $v_d \leq \frac{1}{d}\sum_m \lambda_m^2 \leq K_2 \frac{\lambda_d^2}{d}$, and 
\[
v_d = \frac{1}{d}\sum_{m=1}^d \lambda_m^2 - \Big( \frac{1}{d} \sum_{m=1}^d \lambda_m \Big)^2 \ge \frac{\lambda_d^2}{d} - \frac{K_1^2 \lambda_d^2}{d^2}.
\]
\end{proof}

\begin{corollary}[Upper Bound for the Uniformly Geometric Class]\label{COR_upper_bound_geometric}
Let $A^{(d)} = i \, \mathrm{diag}(\lambda_1, \dots, \lambda_d)$ be a sequence of deterministic drift matrices ordered such that $\vert{}\lambda_1\vert{} \le \dots \le \vert{}\lambda_d\vert{}$. Assume the spectrum belongs to the uniformly geometric class $\mathcal{G}(a,b)$ (Assumption \ref{ASS_geometric_class}).
For any $\epsilon \in (0,1)$, there exist $C_\epsilon>0$ and $d_1 \in \mathbb{N}$ such that, for every $d \geq d_1$,
$$\mathbb{P} \left( \tau_{SS}^{(d)} \le C_\epsilon \frac{d^2}{\vert{}\lambda_d\vert{}} \right) \ge 1 - \epsilon.$$
\end{corollary}

\begin{proof}
By Lemma \ref{LEM_variance_geometric}, $v_d \geq \frac{\lambda_d^2}{2d}$ for $d$ large enough, that is, $\frac{d^{3/2}}{\sqrt{v_d}} \leq \sqrt{2} \frac{d^2}{|\lambda_d|}$. The result follows from Remark \ref{REM_explicit_upper_bound}.
\end{proof}

\begin{proposition}[Spectral Cut Lower Bound for the Uniformly Geometric Class]\label{PRO_quantum_speed_limit_bound_geometric_class}
Assume the spectrum satisfies the Uniformly Geometric Class condition (Assumption \ref{ASS_geometric_class}). Let $\tau_{SS}^{(d)}$ be the state-to-state control time of the randomized system.
For any $\epsilon \in (0,1)$, there exists $d_1 \in \mathbb{N}$ such that, for every $d \geq d_1$,
$$\mathbb{P}\left( \tau_{SS}^{(d)} \ge 2\sqrt{\frac{\epsilon}{K_2}} \frac{\sqrt{d}}{\vert{}\lambda_d\vert{}} \right) \ge 1 - \epsilon,$$
with $K_2$ as in Lemma \ref{LEM_variance_geometric}.
\end{proposition}

\begin{proof}
By Lemma \ref{LEM_variance_geometric}, $\frac{1}{\sqrt{v_d}} \geq \frac{1}{\sqrt{K_2}}\frac{\sqrt{d}}{|\lambda_d|}$ for $d \geq d_1$. The result follows from Proposition \ref{PRO_Uniform_lower_bound}.
\end{proof}

\begin{proof}[Proof of Theorem \ref{THE_random_encadrement_geometric}]
Apply Corollary \ref{COR_upper_bound_geometric} and Proposition \ref{PRO_quantum_speed_limit_bound_geometric_class} with $\epsilon/2$, and use the union bound.
\end{proof}

\begin{remark}\em
The lower bound of Proposition \ref{PRO_quantum_speed_limit_bound_geometric_class} uses a fixed star cut. Optimizing over the vertex $j$ of the star cut (as done in Section \ref{SEC_Time_estimate_dominance} in the dominance regime) is expected to improve it, since the capacity $C_j$ is a weighted sum of finitely many (up to geometrically small terms) almost independent exponential variables.
\end{remark}

\subsection{Time estimates under the spectral dominance assumption}
\label{SEC_Time_estimate_dominance}

We consider now matrices $A^{(d)}$ for which the spectrum satisfies the dominant eigenvalue condition (Assumption \ref{ASS_spectral_dominance}). We denote by $S_{d-1} = \sum_{m=1}^{d-1} \lambda_m^2$ the $\ell^2$ mass of the sub-dominant part of the spectrum, so that $S_{d-1} = o(\lambda_d^2)$. The upper bound of Theorem \ref{THE_random_encadrement_dominance} follows from the general estimate of Theorem \ref{TH_Uniform_bound_MST} once $v_d$ is estimated (Lemma \ref{LEM_variance_dominance}). The lower bound relies on a star cut isolating the vertex $j$ for which $|P_{dj}|$ is minimal (Proposition \ref{PRO_lower_bound_dominance}).

\begin{lemma}\label{LEM_super_exponential_growth}
Let $(\lambda_d)_{d \ge 1}$ be a sequence of non-zero real numbers ordered such that $|\lambda_1| \le |\lambda_2| \le \dots \le |\lambda_d|$. Assume the sequence satisfies the absolute spectral dominance condition:
$$\lim_{d \to \infty} \frac{\lambda_d^2}{\sum_{m=1}^d \lambda_m^2} = 1.$$
Then, the sequence grows super-exponentially. Specifically, for any arbitrary constant $\alpha > 1$, there exists a constant $C > 0$ and an integer $d_0 \ge 1$ such that for all $d \ge d_0$:
$$|\lambda_d| \ge C \cdot \alpha^d.$$
\end{lemma}

\begin{proof}
Let $S_d = \sum_{m=1}^d \lambda_m^2$. The absolute spectral dominance condition is given by:
$$\lim_{d \to \infty} \frac{\lambda_d^2}{S_d} = 1.$$

Using the algebraic identity $\lambda_d^2 = S_d - S_{d-1}$, the limit evaluates to:
\[\lim_{d \to \infty} \left( 1 - \frac{S_{d-1}}{S_d} \right) = 1 \quad \implies \quad \lim_{d \to \infty} \frac{S_{d-1}}{S_d} = 0.\]

By definition, the cumulative sum satisfies $S_{d-1} \ge \lambda_{d-1}^2$. Applying this inequality to the numerator yields:
$$0 \le \frac{\lambda_{d-1}^2}{S_d} \le \frac{S_{d-1}}{S_d}.$$

By the squeeze theorem, we obtain $\lim_{d \to \infty} \frac{\lambda_{d-1}^2}{S_d} = 0$.

Since $\frac{\lambda_{d-1}^2}{\lambda_d^2} = \frac{\lambda_{d-1}^2}{S_d} \frac{S_d}{\lambda_d^2}$ and the second factor converges to $1$ by hypothesis, the product of the limits yields:

$$\lim_{d \to \infty} \frac{\lambda_{d-1}^2}{\lambda_d^2} = 0.$$

Taking the inverse of the square root implies that the quotient of consecutive absolute values diverges to infinity:

$$\lim_{d \to \infty} \frac{|\lambda_d|}{|\lambda_{d-1}|} = +\infty.$$

Let $\alpha > 1$ be an arbitrary constant. The divergence of the sequence of quotients guarantees the existence of an integer $d_0$ such that for all $k \ge d_0$:
$$\frac{|\lambda_k|}{|\lambda_{k-1}|} \ge \alpha \quad \implies \quad |\lambda_k| \ge \alpha |\lambda_{k-1}|.$$

By finite induction on the index $d \ge d_0$, we recursively apply this lower bound:
$$|\lambda_d| \ge \alpha |\lambda_{d-1}| \ge \alpha^2 |\lambda_{d-2}| \ge \dots \ge \alpha^{d-d_0} |\lambda_{d_0}|.$$
By isolating the dimension-dependent exponent from the initial terms, we can write:
$$|\lambda_d| \ge \left( |\lambda_{d_0}| \alpha^{-d_0} \right) \alpha^d.$$
Defining the strictly positive deterministic constant $C = |\lambda_{d_0}| \alpha^{-d_0}$, we establish the geometric lower bound for all $d \ge d_0$:
$|\lambda_d| \ge C \cdot \alpha^d$.
\end{proof}

\begin{lemma}\label{LEM_variance_dominance}
Under Assumption \ref{ASS_spectral_dominance}, $v_d \sim \frac{\lambda_d^2}{d}$ as $d \to \infty$.
\end{lemma}

\begin{proof}
We have $\frac{1}{d}\sum_{m=1}^d \lambda_m^2 = \frac{\lambda_d^2 + S_{d-1}}{d} \sim \frac{\lambda_d^2}{d}$. By the Cauchy-Schwarz inequality, $\sum_{m<d} |\lambda_m| \leq \sqrt{(d-1) S_{d-1}}$, hence
$$\bar{\lambda}_d^2 \leq \frac{1}{d^2}\Big(|\lambda_d| + \sqrt{(d-1)S_{d-1}}\Big)^2 \leq \frac{2\lambda_d^2}{d^2} + \frac{2 S_{d-1}}{d} = o\Big(\frac{\lambda_d^2}{d}\Big).$$
The result follows from $v_d = \frac{1}{d}\sum_m \lambda_m^2 - \bar{\lambda}_d^2$.
\end{proof}

\subsubsection{Asymptotic structure of the couplings}

The following result describes the limit law of a single coupling. It is not needed for the proof of Theorem \ref{THE_random_encadrement_dominance}, but it explains the rank-one structure of the interaction graph in the dominance regime.

\begin{proposition}\label{PRO_extreme_growth_modulus}
Under the absolute spectral dominance assumption (Assumption \ref{ASS_spectral_dominance}), we define the standardized moduli of the off-diagonal elements ($j \neq k$) as:
\[M_{jk}^{(d)} = \frac{d}{|\lambda_d|} |\tilde{A}_{jk}|.\]
Then, as $d \to \infty$, $M_{jk}^{(d)}$ converges in distribution to the product of two independent Rayleigh variables: $M_{jk}^{(d)} \xrightarrow{\mathcal{L}} Z = R_1 R_2$,
where $R_1, R_2$ have the density $f_R(r) = 2r e^{-r^2} \mathbf{1}_{\{r \ge 0\}}$.

The limiting variable $Z$ admits the density:
$$f_Z(z) = 4z K_0(2z) \mathbf{1}_{\{z \ge 0\}},$$
where $K_0$ is the modified Bessel function of the second kind of order 0. Its first moments are $\mathbb{E}[Z] = \frac{\pi}{4}$ and $\mathbb{E}[Z^2] = 1$.
\end{proposition}

\begin{proof}
We isolate the contribution of the dominant eigenvalue by writing $W_{jk}^{(d)} = \frac{d}{\lambda_d} \tilde{A}_{jk}$ as
$$W_{jk}^{(d)} = i d\, \overline{P}_{dj} P_{dk} + R_{jk}^{(d)}, \qquad R_{jk}^{(d)} = i \frac{d}{\lambda_d} \sum_{m=1}^{d-1} \lambda_m \overline{P}_{mj} P_{mk}.$$
Applying \eqref{EQ_variance_general_weights} (Remark \ref{REM_variance_general_weights}) with the weights $\nu = (\lambda_1, \dots, \lambda_{d-1}, 0)$, we get
$$\mathbb{E}\left[|R_{jk}^{(d)}|^2\right] = \frac{d^2}{\lambda_d^2}\, \frac{1}{d^2-1}\Big(S_{d-1} - \frac{1}{d}\Big(\sum_{m=1}^{d-1}\lambda_m\Big)^2\Big) \leq \frac{d^2}{d^2-1}\, \frac{S_{d-1}}{\lambda_d^2} \xrightarrow[d \to \infty]{} 0.$$
The remainder therefore converges to $0$ in mean square, and consequently in probability.

The conjugate of the $d$-th row of $P$ is the $d$-th column of the Haar distributed matrix $P^\dagger$, hence is uniformly distributed on the complex unit sphere. By Theorem \ref{THE_Borel_Poincare} (and invariance under permutation of the coordinates), $(\sqrt{d} P_{dj}, \sqrt{d} P_{dk})$ converges in distribution to a pair $(X, Y)$ of independent $\mathcal{CN}(0,1)$ variables. By Slutsky's theorem, $W_{jk}^{(d)} \xrightarrow{\mathcal{L}} i\, \overline{X} Y$ and $M_{jk}^{(d)} = |W_{jk}^{(d)}| \xrightarrow{\mathcal{L}} |X||Y|$.

Since $X$ and $Y$ are standard complex Gaussians, $R_1 = |X|$ and $R_2 = |Y|$ are independent with density $2r e^{-r^2}$ on $\mathbb{R}_+$. The density of $Z = R_1 R_2$ is
$$f_Z(z) = \int_0^\infty f_R(r) f_R(z/r) \frac{dr}{r} = 4z \int_0^\infty e^{-r^2 - z^2/r^2} \frac{dr}{r} = 2z \int_0^\infty e^{-s - z^2/s} \frac{ds}{s} = 4z K_0(2z).$$
Finally $\mathbb{E}[Z] = \mathbb{E}[R_1]^2 = (\sqrt{\pi}/2)^2 = \pi/4$ and $\mathbb{E}[Z^2] = \mathbb{E}[R_1^2]^2 = 1$.
\end{proof}

\begin{remark}\label{REM_star_heuristic}\em
Proposition \ref{PRO_extreme_growth_modulus} suggests the rank-one approximation $|\tilde{A}_{jk}| \approx \frac{|\lambda_d|}{d} R_j R_k$, with $R_j = \sqrt{d}|P_{dj}|$. Note that the vertices of the interaction graph are indexed by the eigenvectors of $B$, while the index $d$ of $P_{dj}$ refers to the dominant eigenvalue of $A$. Within this approximation, the star spanning tree centered at an arbitrary vertex $k$ has weight 
$\sum_{j \neq k} |\tilde{A}_{jk}|^{-1} \approx \frac{d}{|\lambda_d|} \frac{1}{R_k} \sum_{j\neq k} \frac{1}{R_j}.$
Since $\mathbb{E}[1/R_j] = \sqrt{\pi}$, the law of large numbers suggests a weight of order $\frac{\sqrt{\pi}\, d^2}{|\lambda_d| R_k}$, consistent with Theorem \ref{THE_random_encadrement_dominance}, while the largest edge weight of the star is governed by a Fréchet law of index $2$ and is of order $d^{3/2}/|\lambda_d|$. Choosing the center $k$ maximizing $|P_{dk}|$ (so that $R_k^2 \sim \ln d$) suggests the improved upper bound $\tau_{SS}^{(d)} = \mathcal{O}_P\big(d^2/(|\lambda_d| \sqrt{\ln d})\big)$. 
Making these heuristics rigorous requires a control of the remainder in the rank-one approximation which is uniform with respect to the $d-1$ edges of the star, and of the dependence between the coordinates of a row of $P$; we do not pursue this here.
\end{remark}

\subsubsection{Lower bound}

\begin{proposition}\label{PRO_lower_bound_dominance}
Let $A = i\,\mathrm{diag}(\lambda_1,\dots,\lambda_d)$ with $|\lambda_1| \le \dots \le |\lambda_d|$, $\lambda_d \neq 0$, and let $S_{d-1} = \sum_{m<d}\lambda_m^2$. Let $P$ be Haar distributed and let $j^* \in \{1,\dots,d\}$ be an index such that $|P_{dj^*}| = \min_j |P_{dj}|$. Then, for every $\epsilon \in (0,1)$ and every $d \geq 12/\epsilon$, with probability at least $1-\epsilon$,
\[
C_{j^*} = \sum_{k \neq j^*} |\tilde{A}_{kj^*}|^2 \leq 2 \ln(3/\epsilon) \frac{\lambda_d^2}{d^2} + \frac{3 S_{d-1}}{\epsilon (d-1)},
\]
and therefore $\tau_{SS}^{(d)} \geq 2\Big(2 \ln(3/\epsilon) \frac{\lambda_d^2}{d^2} + \frac{3 S_{d-1}}{\epsilon (d-1)}\Big)^{-1/2}$.

In particular, under the Quantitative Spectral Dominance condition (Assumption \ref{ASS_quantitative_dominance}), for every $\epsilon \in (0,1)$ and $d \geq \max(d_0, 12/\epsilon)$,
\[
\mathbb{P}\left( \tau_{SS}^{(d)} \geq c_\epsilon \frac{d}{|\lambda_d|} \right) \geq 1 - \epsilon, \qquad c_\epsilon = \frac{2}{\sqrt{2\ln(3/\epsilon) + 6K/\epsilon}}.
\]
\end{proposition}

\begin{proof}
As in the proof of Theorem \ref{TH_coupe_minimale_Lyapunov}, $C_{j^*}$ is the variance of the spectrum with respect to the probability weights $(|P_{mj^*}|^2)_m$, hence
$$C_{j^*} \leq \sum_{m=1}^d \lambda_m^2 |P_{mj^*}|^2 = \lambda_d^2 |P_{dj^*}|^2 + \sum_{m=1}^{d-1} \lambda_m^2 |P_{mj^*}|^2.$$
We bound the two terms separately.

\emph{First term.} The conjugate of the $d$-th row of $P$ is uniformly distributed on the complex unit sphere. Hence $(|P_{dj}|^2)_{j} \stackrel{\mathcal{L}}{=} (E_j/\Sigma_E)_j$, where $E_1,\dots,E_d$ are independent $\mathrm{Exp}(1)$ variables and $\Sigma_E = \sum_l E_l$. The minimum $\min_j E_j$ is exponentially distributed with parameter $d$, so that, for $t>0$,
$$\mathbb{P}\Big(\min_j |P_{dj}|^2 > \frac{t}{d^2}\Big) \leq \mathbb{P}\Big(d \min_j E_j > \frac{t}{2}\Big) + \mathbb{P}\Big(\Sigma_E < \frac{d}{2}\Big) \leq e^{-t/2} + \frac{4}{d},$$
by Chebyshev's inequality ($\mathbb{E}[\Sigma_E] = \text{Var}(\Sigma_E) = d$). With $t = 2\ln(3/\epsilon)$ and $d \geq 12/\epsilon$, this probability is at most $2\epsilon/3$.

\emph{Second term.} Let $\mathcal{F}$ be the $\sigma$-algebra generated by the $d$-th row of $P$; the index $j^*$ is $\mathcal{F}$-measurable. For $U \in U(d-1)$, the matrix $\mathrm{diag}(U,1) P$ is Haar distributed and has the same $d$-th row as $P$. Hence, conditionally on $\mathcal{F}$, the law of the vector $q = (P_{mj^*})_{1\le m \le d-1} \in \mathbb{C}^{d-1}$ is invariant under $U(d-1)$. Since $\|q\|^2 = 1 - |P_{dj^*}|^2 \leq 1$, we get $\mathbb{E}[|P_{mj^*}|^2 \mid \mathcal{F}] = \frac{1 - |P_{dj^*}|^2}{d-1} \leq \frac{1}{d-1}$ for $m<d$, and
$$\mathbb{E}\Big[\sum_{m=1}^{d-1} \lambda_m^2 |P_{mj^*}|^2\Big] \leq \frac{S_{d-1}}{d-1}.$$
By Markov's inequality, the second term exceeds $\frac{3S_{d-1}}{\epsilon(d-1)}$ with probability at most $\epsilon/3$.

The first claim follows from the union bound, and the lower bound on $\tau_{SS}^{(d)}$ from Proposition \ref{PRO_rappel_estimations_tauSS} applied with $j = j^*$. Under Assumption \ref{ASS_quantitative_dominance}, $\frac{3S_{d-1}}{\epsilon(d-1)} \leq \frac{6K}{\epsilon}\frac{\lambda_d^2}{d^2}$ since $d-1 \geq d/2$, which gives the second claim.
\end{proof}

\begin{remark}\label{REM_quantitative_dominance_needed}\em
Proposition \ref{PRO_lower_bound_dominance} shows that, in general, $\tau_{SS}^{(d)} \geq c_\epsilon \min\big(\frac{d}{|\lambda_d|}, \sqrt{\frac{d}{S_{d-1}}}\big)$ with probability at least $1-\epsilon$. The second term is irrelevant only if $S_{d-1} = \mathcal{O}(\lambda_d^2/d)$, which is Assumption \ref{ASS_quantitative_dominance}. This restriction is not an artifact of the proof. Consider for instance a spectrum with $|\lambda_{d-1}|/|\lambda_d| = 1/\ln d$ (e.g. $\lambda_d = \prod_{k=2}^d \ln k$), which satisfies Assumption \ref{ASS_spectral_dominance} but not Assumption \ref{ASS_quantitative_dominance}. For the vertex $j^*$ minimizing $|P_{dj}|$, the term $\lambda_{d-1}^2 |P_{d-1,j^*}|^2$, of order $\lambda_{d-1}^2/d$, dominates $\lambda_d^2 |P_{dj^*}|^2$, of order $\lambda_d^2/d^2$. Optimizing over $j$ the sum $\lambda_d^2|P_{dj}|^2 + \lambda_{d-1}^2|P_{d-1,j}|^2$ heuristically leads to a minimal star cut of order $\lambda_d^2 / (d^{3/2}\ln d)$, so that the cut lower bound is of order $d^{3/4}\sqrt{\ln d}/|\lambda_d| \ll d/|\lambda_d|$.
\end{remark}

\begin{proof}[Proof of Theorem \ref{THE_random_encadrement_dominance}]
\emph{Upper bound.} By Lemma \ref{LEM_variance_dominance}, $\frac{d^{3/2}}{\sqrt{v_d}} \leq \sqrt{2}\frac{d^2}{|\lambda_d|}$ for $d$ large enough, and Remark \ref{REM_explicit_upper_bound} gives $\mathbb{P}(\tau_{SS}^{(d)} \le \sqrt{2} C_\epsilon d^2/|\lambda_d|) \geq 1-\epsilon$. By Lemma \ref{LEM_super_exponential_growth} (with $\alpha = 2$), $d^2/|\lambda_d| \to 0$, hence $\tau_{SS}^{(d)} \to 0$ in probability.

\emph{Lower bound.} Under Assumption \ref{ASS_quantitative_dominance}, the lower bound follows from Proposition \ref{PRO_lower_bound_dominance}. The two-sided estimate follows from the union bound, applying each estimate with $\epsilon/2$.
\end{proof}

\section{Asymptotic Classification: Spectral Variance and Localization}
\label{SEC_Classification}

\subsection{Spectral Scale and Relative Concentration}

\begin{proposition}\label{PRO_universal_variance_scaling}
Let $A = i \, \mathrm{diag}(\lambda_1, \dots, \lambda_d)$ be the drift operator with a real spectrum and empirical spectral variance $v_d$, and let $\tilde{A} = P^\dagger A P$, with $P$ Haar distributed. For $j \neq k$,
$$\mathbb{E}[|\tilde{A}_{jk}|^2] = \frac{d}{d^2-1} v_d.$$
Consequently, the characteristic edge control time $\tau_{\mathrm{scale}}(A) = \left( \mathbb{E}[|\tilde{A}_{jk}|^2] \right)^{-1/2}$ satisfies $\tau_{\mathrm{scale}}(A) \sim \sqrt{d/v_d}$ as $d \to \infty$, and, for any $\alpha > 0$,
$\tau_{\mathrm{scale}}(\alpha A) = \frac{1}{\alpha} \tau_{\mathrm{scale}}(A)$.
\end{proposition}

\begin{proof}
The first identity is Lemma \ref{LEM_variance} (which holds for any pair $j \neq k$ by permutation invariance). The scaling $\lambda_m \mapsto \alpha\lambda_m$ maps $v_d$ to $\alpha^2 v_d$.
\end{proof}

All the time estimates of Section \ref{SEC_Time_estimates} can thus be expressed in units of the characteristic scale $\sqrt{d/v_d}$. The remaining dependence on the spectrum is through its \emph{shape}.

\begin{definition}[Shape vector and concentration parameters]
Let $A = i \, \mathrm{diag}(\lambda_1, \dots, \lambda_d)$ be the drift operator with $v_d>0$. 
The shape vector is defined by $\mu_m^{(d)} = \frac{\lambda_m - \bar{\lambda}_d}{\sqrt{d \cdot v_d}}$, $1 \leq m \leq d$. We define the spectral concentration parameter $Q_d$ and the inverse participation ratio $I_d$ of the spectrum by
$$Q_d = \|\boldsymbol{\mu}^{(d)}\|_\infty^2 = \max_m (\mu_m^{(d)})^2, \qquad I_d = \|\boldsymbol{\mu}^{(d)}\|_4^4 = \sum_{m=1}^d (\mu_m^{(d)})^4.$$
\end{definition}

\begin{proposition}\label{PRO_shape_vector}
\begin{enumerate}
\item The shape vector is centered: $\sum_m \mu_m^{(d)} = 0$, and belongs to the unit sphere of $\mathbb{R}^d$: $\|\boldsymbol{\mu}^{(d)}\|_2 = 1$. In particular, $Q_d \in [1/d, 1]$.
\item $Q_d^2 \leq I_d \leq Q_d$. 
\item The Lyapunov condition (Assumption \ref{ASS_Lyapunov_condition}) is equivalent to $I_d \to 0$, hence to $Q_d \to 0$.
\item The Absolute Spectral Dominance condition (Assumption \ref{ASS_spectral_dominance}) implies $Q_d \to 1$.
\end{enumerate}
\end{proposition}

\begin{proof}
(1) By definition of $\bar{\lambda}_d$ and $v_d$, $\sum_m (\lambda_m - \bar{\lambda}_d) = 0$ and $\sum_m (\lambda_m - \bar{\lambda}_d)^2 = d v_d$. Since $\sum_m \mu_m^2 = 1$, $\max_m \mu_m^2 \geq 1/d$.

(2) $\max_m \mu_m^4 \leq \sum_m \mu_m^4 \leq \max_m \mu_m^2 \sum_m \mu_m^2$.

(3) $I_d = \frac{\sum_m (\lambda_m - \bar{\lambda}_d)^4}{d^2 v_d^2} = \frac{\frac{1}{d}\sum_m (\lambda_m - \bar{\lambda}_d)^4}{d v_d^2}$, which is the Lyapunov ratio. The equivalence with $Q_d \to 0$ follows from (2).

(4) By the proof of Lemma \ref{LEM_variance_dominance}, $d v_d \sim \lambda_d^2$ and $\bar{\lambda}_d^2 = o(\lambda_d^2)$. Hence $(\mu_d^{(d)})^2 = \frac{(\lambda_d - \bar{\lambda}_d)^2}{d v_d} \to 1$, and $Q_d \to 1$.
\end{proof}

\begin{proposition}\label{PRO_universal_factorization}
Let $A = i \, \mathrm{diag}(\lambda_1, \dots, \lambda_d)$ be the drift operator, and let $\tilde{A}_{jk} = \langle P_{\cdot, j}, A P_{\cdot, k} \rangle$ be the off-diagonal elements of the interaction matrix rotated by a Haar-distributed unitary matrix $P$.

The inverse edge weights $W_{jk}^{(d)} = |\tilde{A}_{jk}|^{-1}$ of the interaction graph admit the  factorization:
$$W_{jk}^{(d)} = \sqrt{\frac{d}{v_d}} |V_{jk}^{(d)}|^{-1},$$
where $V_{jk}^{(d)}$ is a standardized complex random variable, with $\mathbb{E}[V_{jk}^{(d)}]=0$ and $\mathbb{E}[|V_{jk}^{(d)}|^2] = \frac{d^2}{d^2-1}$ (Lemma \ref{LEM_variance}), defined by:
$$V_{jk}^{(d)} = i \sum_{m=1}^d \mu_m^{(d)} \left( \sqrt{d} P_{mj}^* \right) \left( \sqrt{d} P_{mk} \right),$$
and $\boldsymbol{\mu}^{(d)} = (\mu_1^{(d)}, \dots, \mu_d^{(d)})$ is the spectral shape vector defined by $\mu_m^{(d)} = \frac{\lambda_m - \bar{\lambda}_d}{\sqrt{d \cdot v_d}}$.
\end{proposition}

\begin{proof}[Proof of Proposition \ref{PRO_universal_factorization}
]
By definition, the rotated interaction matrix elements for $j \neq k$ are given by:

$$\tilde{A}_{jk} = i \sum_{m=1}^d \lambda_m P_{mj}^* P_{mk}.$$

Because the columns of the unitary matrix $P$ are orthogonal, 
$$\langle P_{\cdot, j}, P_{\cdot, k} \rangle = \sum_{m=1}^d P_{mj}^* P_{mk} = 0.$$ 
Consequently, shifting the entire spectrum by any global scalar $c \in \mathbb{R}$ leaves the off-diagonal elements invariant:

$$\sum_{m=1}^d (\lambda_m - c) P_{mj}^* P_{mk} = \sum_{m=1}^d \lambda_m P_{mj}^* P_{mk} - c \langle P_{\cdot, j}, P_{\cdot, k} \rangle = \sum_{m=1}^d \lambda_m P_{mj}^* P_{mk}.$$

We select $c = \bar{\lambda}_d$ to center the spectrum. The element becomes:

$$\tilde{A}_{jk} = i \sum_{m=1}^d (\lambda_m - \bar{\lambda}_d) P_{mj}^* P_{mk}.$$

We explicitly factor out the macroscopic scaling term $\sqrt{v_d / d}$ by multiplying and dividing the summand by $d \sqrt{v_d}$. We also extract the dimension scaling $1/d$ to normalize the eigenvector components:

$$\tilde{A}_{jk} = i \sqrt{\frac{v_d}{d}} \sum_{m=1}^d \left( \frac{\lambda_m - \bar{\lambda}_d}{\sqrt{d \cdot v_d}} \right) \left( \sqrt{d} P_{mj} \right)^* \left( \sqrt{d} P_{mk} \right).$$

By substituting the definition of the spectral shape vector components $\mu_m^{(d)} = \frac{\lambda_m - \bar{\lambda}_d}{\sqrt{d \cdot v_d}}$, we identify the standardized variable:

$$\tilde{A}_{jk} = \sqrt{\frac{v_d}{d}} V_{jk}^{(d)} \quad \text{where} \quad V_{jk}^{(d)} = i \sum_{m=1}^d \mu_m^{(d)} X_{mj}^* X_{mk},$$

and $X_{mj} = \sqrt{d} P_{mj}$. Taking the inverse absolute value yields the factorization $W_{jk}^{(d)} = \sqrt{\frac{d}{v_d}} |V_{jk}^{(d)}|^{-1}$. Note that $V_{jk}^{(d)} = Z_{j,k}^{(d)}$ with the notation of Proposition \ref{PRO_off_diagonal_normality}.

\end{proof}

\begin{proposition}\label{PRO_asymptotic_behavior_universality_parameter}
Let $j \neq k$.
\begin{enumerate}
\item If $Q_d \to 0$ (Lyapunov regime), then $V_{jk}^{(d)} \xrightarrow{\mathcal{L}} \mathcal{CN}(0,1)$.
\item If $Q_d \to 1$ (dominance regime), then $V_{jk}^{(d)} \xrightarrow{\mathcal{L}} i\,\overline{X} Y$, where $X,Y$ are independent $\mathcal{CN}(0,1)$ variables. In particular $|V_{jk}^{(d)}| \xrightarrow{\mathcal{L}} R_1 R_2$, the product of two independent Rayleigh variables of Proposition \ref{PRO_extreme_growth_modulus}.
\end{enumerate}
\end{proposition}

\begin{proof}
(1) By Proposition \ref{PRO_shape_vector}, $Q_d \to 0$ is equivalent to the Lyapunov condition, and $V_{jk}^{(d)} = Z_{j,k}^{(d)}$. The result is Proposition \ref{PRO_off_diagonal_normality}.

(2) For every $d$, let $m^*_d$ be an index such that $(\mu^{(d)}_{m^*_d})^2 = Q_d$, and write
$$V_{jk}^{(d)} = i\, \mu^{(d)}_{m^*_d}\, \overline{X}_{m^*_d j} X_{m^*_d k} + R^{(d)}, \qquad R^{(d)} = i \sum_{m \neq m^*_d} \mu^{(d)}_m\, \overline{X}_{mj} X_{mk}.$$
By \eqref{EQ_variance_general_weights} (Remark \ref{REM_variance_general_weights}), applied with the weights $\nu_m = \mu_m^{(d)}$ for $m \neq m^*_d$ and $\nu_{m^*_d} = 0$,
$$\mathbb{E}[|R^{(d)}|^2] \leq \frac{d^2}{d^2-1} \sum_{m \neq m^*_d} (\mu^{(d)}_m)^2 = \frac{d^2}{d^2-1}(1 - Q_d) \xrightarrow[d\to\infty]{} 0.$$
For the main term, let $\sigma_d = \mathrm{sign}(\mu^{(d)}_{m^*_d}) \in \{-1,1\}$. Permuting the rows of $P$ (so as to bring the $m^*_d$-th row to the first position) and multiplying its $j$-th column by $\sigma_d$ (operations which preserve the Haar measure), we see that $\sigma_d \overline{X}_{m^*_d j} X_{m^*_d k}$ has the same law as $\overline{X}_{1j} X_{1k}$. By Theorem \ref{THE_Borel_Poincare} applied to the first column of the Haar matrix $P^\dagger$ (and permutation invariance), $(\sqrt{d}P_{1j}, \sqrt{d}P_{1k}) \xrightarrow{\mathcal{L}} (X, Y)$, so that $\overline{X}_{1j} X_{1k} \xrightarrow{\mathcal{L}} \overline{X}Y$. Since $|\mu^{(d)}_{m^*_d}| = \sqrt{Q_d} \to 1$, Slutsky's theorem gives the result.
\end{proof}

\begin{remark}\em
The summands $\overline{X}_{mj} X_{mk}$, $m = 1,\dots,d$, of $V_{jk}^{(d)}$ are not independent (Remark \ref{REM_variance_general_weights}): they sum to $d\langle P_{\cdot,j},P_{\cdot,k}\rangle = 0$, and their pairwise covariances, of order $1/d$, have a cumulative effect of order $1$ on non-centered weighted sums. For this reason, Proposition \ref{PRO_asymptotic_behavior_universality_parameter} cannot be obtained from the Lindeberg-Feller central limit theorem for independent variables. In case (1), the conditioning argument of Section \ref{SEC_Lyapunov_normality} takes the orthogonality of the columns into account; in case (2), only the second moment \eqref{EQ_variance_general_weights} is needed. In both cases, the fact that the shape vector is centered is essential (consider the non-centered weights $\nu_m = 1/\sqrt{d}$, for which $\sum_m \nu_m \overline{X}_{mj}X_{mk} = 0$ identically).
\end{remark}

In the intermediate uniformly geometric class, $Q_d$ stays away from $0$ and $1$ (Proposition \ref{PRO_IPR_geometric_class} below), and the limit laws of $V_{jk}^{(d)}$ are infinite weighted sums of products of independent complex Gaussian variables. Since the spectral ratios are only assumed to lie in $[a,b]$, the limit is defined only along subsequences.

\begin{proposition}\label{PRO_subsequential_convergence}
Assume that the spectrum belongs to the uniformly geometric class $\mathcal{G}(a,b)$ and let $j \neq k$. For $0 \le m \le d-1$, let $c_m^{(d)} = \mu_{d-m}^{(d)}$.
For any subsequence $(d_l)_{l \ge 1}$, there exist a further subsequence $(d_{l_r})_{r \ge 1}$ and a real sequence $\boldsymbol{\gamma} = (\gamma_m)_{m \ge 0}$ such that $c_m^{(d_{l_r})} \to \gamma_m$ for every $m \geq 0$,
$$\sum_{m \ge 0} \gamma_m^2 = 1, \qquad \sqrt{1-a^{-2}}\; b^{-m} \le |\gamma_m| \le \sqrt{1-b^{-2}}\; a^{-m} \quad \forall m \ge 0,$$
and
$$V_{jk}^{(d_{l_r})} \xrightarrow{\mathcal{L}} Z_{\boldsymbol{\gamma}} = i \sum_{m=0}^\infty \gamma_m \overline{X}_m Y_m,$$
where $(X_m)_{m\ge0}$, $(Y_m)_{m\ge0}$ are independent $\mathcal{CN}(0,1)$ variables (the series converges in $L^2$).
\end{proposition}

\begin{proof}
We first establish the bounds on the coefficients. By Lemma \ref{LEM_variance_geometric}, $\bar{\lambda}_d^2 \le K_1^2 \lambda_d^2/d^2$ and $d v_d \geq \lambda_d^2 (1 - K_1^2/d)$, so that $\bar{\lambda}_d^2/(d v_d) = \mathcal{O}(1/d^2)$ and $d \bar{\lambda}_d^2/(d v_d) = \mathcal{O}(1/d)$. Moreover, $d v_d = \sum_m \lambda_m^2 - d\bar{\lambda}_d^2 = \sum_m \lambda_m^2 + \mathcal{O}(\lambda_d^2/d)$. Hence, for fixed $m$ and $d \to \infty$,
$$c_m^{(d)} = \frac{\lambda_{d-m}}{\sqrt{\sum_{n} \lambda_n^2}} + o(1), \qquad |c_m^{(d)}| = \frac{|\lambda_{d-m}|}{|\lambda_d|} \sqrt{Q'_d} + o(1), \qquad \mbox{where } Q'_d = \frac{\lambda_d^2}{\sum_n \lambda_n^2}.$$
For $d - m \geq m_0$, we have $b^{-m} \le |\lambda_{d-m}|/|\lambda_d| \le a^{-m}$, and the proof of Proposition \ref{PRO_IPR_geometric_class} shows that all the accumulation points of $Q'_d$ lie in $[1-a^{-2}, 1-b^{-2}]$. Since each sequence $(c_m^{(d)})_d$ is bounded, Cantor's diagonal argument provides a subsequence along which $c_m^{(d)} \to \gamma_m$ for every $m$ and $Q'_d$ converges, and the limits satisfy the announced bounds.

We now control the tails of the coefficient sequences uniformly in $d$. For $M \ge 0$, using the inequality $(x-y)^2 \le 2x^2 + 2y^2$, we get
$$\sum_{M < m \le d-1} (c_m^{(d)})^2 \le \frac{2}{d v_d} \sum_{n < d-M} \lambda_n^2 + \frac{2 d \bar{\lambda}_d^2}{d v_d} \le C a^{-2M} + \frac{C}{d},$$
for a constant $C$ independent of $d$ and $M$, by the argument of Lemma \ref{LEM_variance_geometric}. Since $\sum_{m} (c_m^{(d)})^2 = 1$ for every $d$, Fatou's lemma gives $\sum_m \gamma_m^2 \le 1$, while the tail estimate gives $\sum_{m \le M} \gamma_m^2 = \lim \sum_{m \le M} (c_m^{(d)})^2 \geq 1 - C a^{-2M}$ along the extracted subsequence, for every $M$. Hence $\sum_m \gamma_m^2 = 1$.

It remains to identify the limit law. Let $X_{m,d} = \sqrt{d}P_{d-m,j}$ and $Y_{m,d} = \sqrt{d}P_{d-m,k}$. For any fixed $M$, Lemma \ref{LEM_block_Gaussian} shows that $(X_{m,d}, Y_{m,d})_{0 \le m \le M}$ converges in distribution to a family $(X_m, Y_m)_{0 \le m \le M}$ of independent $\mathcal{CN}(0,1)$ variables. Along the extracted subsequence, the continuous mapping theorem then gives
$$V^{(d)}_{\le M} := i \sum_{m=0}^M c_m^{(d)} \overline{X}_{m,d} Y_{m,d} \xrightarrow{\mathcal{L}} i \sum_{m=0}^M \gamma_m \overline{X}_m Y_m.$$
On the one hand, by \eqref{EQ_variance_general_weights} and the tail estimate above, $\mathbb{E}|V_{jk}^{(d)} - V^{(d)}_{\le M}|^2 \le \frac{d^2}{d^2-1}\sum_{m>M} (c_m^{(d)})^2 \le 2C(a^{-2M} + 1/d)$, so that, by Chebyshev's inequality, $\lim_{M} \limsup_{d} \mathbb{P}(|V_{jk}^{(d)} - V^{(d)}_{\le M}| > \eta) = 0$ for every $\eta > 0$. On the other hand, the variables $\overline{X}_m Y_m$ being independent, centered and of unit variance, $\mathbb{E}|Z_{\boldsymbol{\gamma}} - i\sum_{m\le M} \gamma_m \overline{X}_m Y_m|^2 = \sum_{m > M}\gamma_m^2 \to 0$ as $M \to \infty$. The standard approximation argument for convergence in distribution (if $U_{M,d} \to U_M$ in law as $d\to\infty$ for each $M$, $U_M \to U$ in law, and $\lim_M \limsup_d \mathbb{P}(|U_{M,d} - U_d| > \eta) = 0$ for every $\eta>0$, then $U_d \to U$ in law; see e.g. \cite{billingsley1995probability}) concludes the proof.
\end{proof}

\begin{proposition}\label{PRO_IPR_geometric_class}
Let $(\lambda_m)_{m \ge 1}$ be the spectrum of the drift operator and $Q_d = \frac{\max_m (\lambda_m - \bar{\lambda}_d)^2}{d \cdot v_d}$ be its concentration parameter at dimension $d$.

\textbf{Pure Geometric Spectrum.}
If the spectrum is purely geometric, $\lambda_m = c \cdot r^m$ with $c \neq 0$ and $r > 1$, then $Q_d$ converges to a constant:
$$\lim_{d \to \infty} Q_d = 1 - r^{-2} \in (0, 1).$$

\textbf{Uniformly Geometric Class.}
If the spectrum belongs to the uniformly geometric class $\mathcal{G}(a,b)$ with $1 < a \le b < \infty$ (as per Assumption \ref{ASS_geometric_class}), then the sequence $Q_d$ is bounded away from 0 and 1 for sufficiently large $d$. Specifically, all its accumulation points $Q^*$ satisfy:
$$1 - a^{-2} \le Q^* \le 1 - b^{-2}.$$
\end{proposition}

\begin{proof}
By Lemma \ref{LEM_variance_geometric} (which applies to the pure geometric case with $a=b=r$, $m_0=1$), $\bar{\lambda}_d = \mathcal{O}(\lambda_d/d)$ and $d v_d = S_2^{(d)} + \mathcal{O}(\lambda_d^2/d)$, where $S_2^{(d)} = \sum_{m=1}^d \lambda_m^2 \geq \lambda_d^2$. The maximal squared deviation is attained at $m=d$ for $d$ large and satisfies $(\lambda_d - \bar{\lambda}_d)^2 \sim \lambda_d^2$. Hence
$$Q_d \sim \frac{\lambda_d^2}{S_2^{(d)}}.$$

For the pure geometric spectrum, $S_2^{(d)} = c^2 r^2 \frac{r^{2d}-1}{r^2-1} \sim \lambda_d^2 \frac{r^2}{r^2-1}$, so that $Q_d \to 1 - r^{-2}$.

For the uniformly geometric class, for $0 \le k \le d - m_0$, we have
$b^{-k} \le \frac{|\lambda_{d-k}|}{|\lambda_d|} \le a^{-k}$, while the terms with $m < m_0$ contribute at most $m_0 a^{-2(d-m_0)} \lambda_d^2$ to $S_2^{(d)}$. Hence
$$\lambda_d^2 \sum_{k=0}^{d-m_0} b^{-2k} \le S_2^{(d)} \le \lambda_d^2 \Big( \sum_{k=0}^{\infty} a^{-2k} + m_0 a^{-2(d-m_0)}\Big),$$
and, letting $d \to \infty$,
$$\frac{1}{1 - b^{-2}} \le \liminf_{d\to\infty} \frac{S_2^{(d)}}{\lambda_d^2} \le \limsup_{d\to\infty} \frac{S_2^{(d)}}{\lambda_d^2} \le \frac{1}{1 - a^{-2}}.$$
Inverting these bounds, we obtain
$$1 - a^{-2} \le \liminf_{d \to \infty} Q_d \le \limsup_{d \to \infty} Q_d \le 1 - b^{-2}.$$
Since $a \le b$, this interval is a nonempty subset of $(0,1)$.
\end{proof}

\begin{remark}
The level of the concentration parameter $Q_d$ is not, by itself, a certificate that the spectrum belongs to one of the three regimes above. Consider the power spectrum $\lambda_k = k^a$ with $a < -1/2$ (Section \ref{SEC_polynomial_spectrum}), reindexed so that $|\lambda_1| \le \dots \le |\lambda_d|$. Since the series $\sum_k k^{2a}$ and $\sum_k k^{4a}$ converge,

$$v_d \sim \frac{\zeta(-2a)}{d}, \qquad \frac{\lambda_d^2}{\sum_{k=1}^d \lambda_k^2} \to \frac{1}{\zeta(-2a)} \in (0,1), \qquad \frac{\tfrac{1}{d}\sum_{k=1}^d \lambda_k^4}{d\, v_d^2} \to \frac{\zeta(-4a)}{\zeta(-2a)^2} \neq 0.$$

The Lyapunov condition fails, the largest eigenvalue does not dominate, and the spectrum is not in the uniformly geometric class either, the infimum of the ratios of consecutive eigenvalues tending to $1$. Yet the concentration parameter converges to the intermediate value $Q_d \to \zeta(-2a)^{-1} \in (0,1)$, at the same level as for the geometric class, and the shape vector (read from the largest eigenvalue downwards) converges in $\ell^2$ to the explicit sequence $\boldsymbol{\gamma}_a = (k^a/\sqrt{\zeta(-2a)})_{k \geq 1} \in \ell^2$.

The same argument as the proof of Proposition \ref{PRO_subsequential_convergence} (here the $\ell^2$ convergence of the shape vector replaces the subsequence extraction) then gives, for $j \neq k$,

$$V_{jk}^{(d)} \xrightarrow{\mathcal{L}} Z_{\boldsymbol{\gamma}_a} := i \sum_{k=1}^{\infty} \gamma_{a,k}\, \overline{X}_k Y_k,$$

an infinite weighted sum of products of independent complex Gaussian variables, of the same type as the limit $Z_{\boldsymbol{\gamma}}$ of Proposition \ref{PRO_subsequential_convergence} but with algebraically, rather than geometrically, decaying coefficients.

Two observations follow. First, $Q_d$ is the right coarse invariant to separate the three regimes, $Q_d \to 0$ (Lyapunov), $Q_d$ bounded away from $0$ and $1$ (geometric class), $Q_d \to 1$ (dominance); the preceding example shows that the intermediate level is not, however, exclusive to the geometric class. Second, $Q_d$ is not fine enough to determine the asymptotic law of the edges: two spectra may satisfy the same relative scaling $v_d \asymp \lambda_d^2/d$ and have the same limit of $Q_d$ and yet have different shapes. For instance, the harmonic spectrum $\lambda_k = 1/k$ (the case $a=-1$) satisfies $Q_d \to 6/\pi^{2}$, and so does the geometric spectrum with ratio $r$ chosen such that $1-r^{-2} = 6/\pi^{2}$ ($r \simeq 1.6$); the corresponding shape vectors have algebraic decay ($k^{-1}$) and geometric decay ($r^{-k}$) respectively, and the limit laws $Z_{\boldsymbol{\gamma}_a}$ and $Z_{\boldsymbol{\gamma}}$ are distinct. The full shape vector $\boldsymbol{\mu}^{(d)}$ is thus the finer invariant on which the control-time estimates ultimately depend.
\end{remark}

\subsection{Concluding Discussion}

\paragraph{Summary of the scalings.}
By Proposition \ref{PRO_universal_variance_scaling}, the natural time unit is the characteristic edge time $\sqrt{d/v_d}$. In this unit, the main results of the paper read as follows (with high probability, up to constants):

\begin{center}
{\small
\begin{tabular}{lccc}
\toprule
Regime & $Q_d$ & lower bound & upper bound \\
\midrule
Lyapunov (Theorem \ref{THE_random_encadrement}) & $\to 0$ & $d^{-1/2}$ & $d$ \\
Uniformly geometric (Theorem \ref{THE_random_encadrement_geometric}) & $\begin{gathered} Q_d \to c \ \ \text{(up to subsequence)}\\ c \in (0,1) \end{gathered}$ & $d^{-1/2}$ & $d$ \\
Quantitative dominance (Theorem \ref{THE_random_encadrement_dominance}) & $\to 1$ & $1$ & $d$ \\
\bottomrule
\end{tabular}
}
\end{center}

The ratio between the upper and the lower bound is therefore of order $d^{3/2}$ in the Lyapunov and uniformly geometric regimes, and of order $d$ in the dominance regime. In the geometric case, every accumulation point of $Q_d$ lies in $[1-a^{-2}, 1-b^{-2}]$ (Proposition \ref{PRO_IPR_geometric_class}); the time bounds hold along the whole sequence, the extraction being only needed for the limit laws.

\paragraph{Upper bounds.} In all regimes, the upper bound is obtained from a fixed spanning tree (Theorem \ref{TH_Uniform_bound_MST}), and the factor $d$ simply counts its edges, each of them having a typical weight of order $\sqrt{d/v_d}$. This bound does not exploit the possibility of selecting strong couplings. In the Lyapunov regime, the heuristics discussed in Section \ref{SEC_polynomial_spectrum} suggest that the minimal spanning tree is lighter by a factor $(\ln d)^{-1/2}$; in the dominance regime, Remark \ref{REM_star_heuristic} suggests the same logarithmic gain for a well chosen star. Whether the true control time is closer to the lower or to the upper bound is left open: the spanning tree upper bound relies on sequential transfers along the edges of the tree (one-dimensional effective transport), while the cut lower bound only accounts for the total capacity of the couplings of a single vertex.

\paragraph{Lower bounds.} The difference between the regimes appears in the lower bound. In the Lyapunov and geometric regimes, the capacity of the star cut centered at a given vertex is of order $v_d$ (in probability in the Lyapunov case, in mean in the geometric case), and the lower bound $\sqrt{d/v_d} \cdot d^{-1/2} = 1/\sqrt{v_d}$ follows. In the dominance regime, the interaction matrix is asymptotically of rank one, $|\tilde{A}_{jk}| \approx |\lambda_d| |P_{dj}||P_{dk}|$, so that the coupling of a vertex $j$ to the rest of the network is entirely controlled by the single coordinate $|P_{dj}|$. The vertex minimizing $|P_{dj}|$ is weakly coupled to all other vertices, which yields the larger lower bound $\sqrt{d/v_d} \asymp d/|\lambda_d|$ (Proposition \ref{PRO_lower_bound_dominance}). This mechanism requires the remainder of the spectrum to be small enough (Assumption \ref{ASS_quantitative_dominance} and Remark \ref{REM_quantitative_dominance_needed}).

\paragraph{Physical interpretation.} Expressed in terms of the largest eigenvalue, the upper bound $d^{3/2}/\sqrt{v_d}$ tends to zero as soon as the spectrum grows fast enough ($\lambda_k = k^a$ with $a>3/2$, geometric or super-exponential spectra). In these cases, increasing the dimension of the Hilbert space decreases the control time with high probability, despite the growing number of states to be reached by a single control.

\section*{Acknowledgments}

The authors gratefully acknowledge the financial support that made this work possible.  
It was funded, in whole or in part, by the French National Research Agency (ANR) under the QuBiCCS project (grant ANR-24-CE40-3008-01).  The research is also part of the EIPHI Graduate School initiative (grant ANR-17-EURE-0002) and received additional backing from the Bourgogne Franche-Comté Region.

The authors also acknowledge the use of Gemini and the large-language-model services
provided by the University of Strasbourg (Unistra) for drafting the manuscript and converting
the original French notes into academic English.  While these tools contributed to the writing, the authors
take full responsibility for the scientific integrity, mathematical correctness, and
final verification of the article.

\bibliographystyle{siam}
\bibliography{Biblio}
\end{document}